\documentclass [11pt]{amsart}
\usepackage {amsmath, amssymb, a4wide, epsfig,enumerate}
\usepackage[latin1]{inputenc}
\date{\today}

\author{Tien Cuong Dinh, Romain Dujardin and Nessim Sibony}
\title{On the dynamics near infinity of some polynomial mappings in $\cd$}

\newcommand{\cc}{\mathbb{C}}
\newcommand{\bb}{\mathbb{B}}

\newcommand{\dd}{\mathbb{D}}

\newcommand{\nn}{\mathbb{N}}
\newcommand{\e}{\varepsilon}
\newcommand{\cv}{\rightarrow}
\newcommand{\cvf}{\rightharpoonup}
\newcommand{\fr}{\partial}

\newcommand{\set}[1]{\left\{#1\right\}}
\newcommand{\norm}[1]{\left\Vert#1\right\Vert}
\newcommand{\abs}[1]{\left\vert#1\right\vert}
\newcommand{\finc}{\subset \subset}
\newcommand{\cd}{\cc^2}
\newcommand{\pd}{{\mathbb{P}^2}}
\newcommand{\pu}{{\mathbb{P}^1}}
\newcommand{\rest}[1]{ \arrowvert_{#1}}

\newcommand{\unsur}[1]{\frac{1}{#1}}
\newcommand{\el}{\mathcal{L}}

\newcommand{\tel}{{}^t \mathcal{L}}

\newcommand{\rond}{\hspace{-.1em}\circ\hspace{-.1em}}

\DeclareMathOperator{\supp}{Supp}

\DeclareMathOperator{\dist}{dist}

\newtheorem{prop} {Proposition} [section]
\newtheorem{thm}[prop] {Theorem} 
\newtheorem{defi}[prop] {Definition}
\newtheorem{lem}[prop] {Lemma}
\newtheorem{cor}[prop]{Corollary}
\newtheorem{theo} {Theorem}
\newtheorem{exam}[prop]{Example}
\newenvironment{example}{\begin{exam}\normalfont}{\end{exam}}
\newtheorem{rmk}[prop]{Remark}

\newenvironment{pf}{\noindent {\bf Proof:}}{\hfill $\square$\\}


\begin{document}

\begin{abstract}
We construct the Green current for  a random iteration of horizontal-like
mappings in $\cd$. This is applied to the study of a polynomial map $f:\cd\cv\cd$
with the following properties:
\begin{itemize}
\item[\em i.] infinity is $f$-attracting;
\item[\em ii.] $f$ contracts the line at infinity to a point not in the indeterminacy set.
\end{itemize}
Then the Green current of $f$ can be decomposed into pieces 
associated with an itinerary determined by the indeterminacy points. We also 
study the escape rates near infinity, i.e. the possible values of the function 
$\limsup \unsur{n}\log^+\log^+ \norm{f^n}$. 
\end{abstract}

\maketitle

\section{Introduction}

Exploring the behavior of a polynomial mappings  $\cd\cv\cd$ as a
dynamical system is an extremely rich topic. Much attention has been paid in the last
few years in constructing currents and measures with interesting dynamical
properties. Parallel to this is the need to find global models for the dynamics.
In the present paper  we study
a ``model at infinity" for certain polynomial maps in $\cd$, with indeterminacy points
at infinity. Our results illustrate the kind of dynamical complexity that can arise
in this setting. On the other hand our method partially relies on
delicate convergence
results for some sequences of positive closed currents, obtained through
iteration along a sequence
 of horizontal-like maps.\\

Let us be more specific. We consider polynomial mappings in 
$\cd$, with meromorphic compactification on $\pd$.
The  indeterminacy set $I(f)=\set{I_1,
\ldots , I_m}$ ($m\geq 2$) is a finite number of points in the line at infinity 
$L_\infty$. We assume $L_\infty$ 
 is $f$-attracting , and $f(L_\infty\setminus I(f))$
is a point
$X\notin I(f)$. Thus for each $j$, $f(I_j)=L_\infty$.
Each indeterminacy point $I_i$ comes equipped with a multiplicity
$d_i$,  which is the intersection multiplicity at $I_i$ of $L_\infty$
and $f^{-1}(L)$, where $L$ is a generic line.

We denote by $T$ the Green current of $f$, which was shown to
exist in this case by the third author \cite{sib}.

For a point $p$ escaping to infinity, we define its {\sl escape rate}
 to be $\limsup \unsur{n}\log^+\log^+ \norm{f^n(p)}$.
We make the
assumption that a local escape rate $\ell_i$ is associated
to  each indeterminacy point, that is, for every $i$ there exists
constants   $C_1, C_2$,
  ${\ell}_i\geq 1$, and neighborhoods $V(I_i)$ and $V(X)$
such that
$$
\forall p\in V(I_i), \text{ s.t. } f(p)\notin V(X),~C_1\norm{p}^{{\ell}_i}
\leq \norm{f(p)} \leq C_2\norm{p}^{{\ell}_i}.
$$

In \cite{dis2} an algebraically generic class of polynomial mappings
$\cc^k\cv\cc^k$ was constructed, and for that class there are only
 finitely many escape
rates to infinity. This was also the case in all previously known
examples. In contrast with these phenomena, our main result is the
following. This is a 
combination of theorem \ref{thm_speed} and
corollary \ref{cor_decompos}.

\begin{theo} Assume  the preceding hypotheses hold. Then
\begin{enumerate}[\em i.]
\item The set of possible escape rates to infinity is
$[\min \ell_i,\max \ell_i]\cup\set{d}$.
\item For $\norm{T}$-a.e. point in the basin of infinity, the limit
$\ell(p)= \lim \unsur{n}\log^+\log^+ \norm{f^n(p)}$ exists and
$\displaystyle{\ell(p)= \prod_{i=1}^m \ell_i^\frac{d_i}{d}}$.
\end{enumerate}
\end{theo}

The key point is that in a  small bidisk neighborhood
$\bb_i\ni I_i$, $f$ acts as a {\sl horizontal-like map} \cite{hl}, that is,
a meromorphic map stretching the horizontal direction and
contracting the vertical
direction in $\bb_i$ --here the line at infinity is viewed
 as horizontal. Such maps act
by pull back on positive closed   currents with vertical support.
 Each point escaping to infinity,
not in the basin of the superattracting point $X$, has an
{\sl itinerary} in the set of boxes $\bb_i$, hence
  in  the symbol space $\Sigma=\set{1, \ldots, m}^\nn$.
We are thus led to consider the iterates of a
vertical current through this  prescribed itinerary.

Proving convergence results for such sequences of currents is a main technical
issue in the paper. We think
this formalism should find numerous applications, some of which are sketched in
\S \ref{sec_appl}. This occurs for instance
in studying the local dynamics around a periodic  indeterminacy point, as shown
in \cite{ya}.

 The Green current then
admits, in the basin of infinity, a decomposition
into pieces that reflect the partition
induced by the itinerary map (theorem \ref{thm_decompos}).
On an algebraically generic class of examples,
the pieces are just
holomorphic disks near infinity (example \ref{ex_externalrays}). This allows
to reconstruct the maximal entropy measure through landing of a family of
external rays (theorem \ref{thm_externalrays}). \\

The precise outline of the paper is as follows. In \S \ref{sec_hl} we define
horizontal-like maps and detail their action on currents. We state the main
convergence
theorem \ref{thm_cv_main} which will be proved in \S \ref{sec_conv}. Section
\ref{sec_speed} is devoted to studying our polynomial maps in $\cd$, whereas
\S \ref{sec_appl} lists some other results and problems involving sequences of
horizontal-like maps.


\section{Meromorphic horizontal-like maps}\label{sec_hl}

We introduce the class of local maps 
that will serve as a building block for the following. Horizontal-like 
maps were defined in \cite{hl}, however the focus in that paper was rather
on biholomorphic maps, whereas we are led here
to consider meromorphic maps. We 
define and study the basic properties of the
Graph Transform operators for currents in this setting.\\

Let us first fix some notation. We let $\bb$ be the unit bidisk in $\cd$, 
and respectively $\fr_v\bb$ and $\fr_h\bb$ the vertical and horizontal part
of the boundary, that is
$$\fr_v\bb=\set{(z,w)\in\cd,~ \abs{z}=1,~\abs{w}<1} \text{ and }
\fr_h\bb=\set{(z,w)\in\cd,~ \abs{z}<1,~\abs{w}=1}.$$
We also denote by $\dd$ (resp. $\dd_r$) the unit disk in $\cc$
(resp. the disk $D(0,r)\subset \cc$).

\begin{defi}\label{def_hl} 
Let $\bb_i\subset M_i$ be an open subset
biholomorphic  to $\bb$ in the complex manifold $M_i$, for $i=1,2$.
Let $f$ be a dominating meromorphic map 
defined in some neighborhood of  $\overline{\bb_1}$ with
values in $M_2$.
 The triple $(f,\bb_1,\bb_2)$ defines a horizontal-like map
 if
 \begin{itemize}
\item[i.] $f$ has no indeterminacy points in $\fr_v\bb_1$ and 
$f(\fr_v\bb_1)\cap \overline{\bb_2}=\emptyset$;
\item[ii.] $f(\overline{\bb_1})\cap\fr\bb_2\subset \fr_v\bb_2$;
\item[iii.] $f(\bb_1)\cap\bb_2\neq \emptyset$.
\end{itemize}
\end{defi}

Some comments are in order here. First, the definition is analogous to the one in
\cite{hl}, except that indeterminacy points are allowed, and the source and target 
bidisks may be different. 
That definition was in turn inspired by the {\sl 
crossed mappings} in \cite{ho}.
We often abuse notation and speak of a horizontal-like map 
$f:\bb_1\cv\bb_2$ even if 
the range of $f$  is not contained  in $\bb_2$. The definition can equally be written
without reference to the ambient manifolds $M_i$ by specifying 
the  properties  of the graph of $f$ in $\bb_1\times \bb_2$. More specifically, 
a 2-dimensional subvariety
$\Gamma\subset\bb_1\times \bb_2$ defines a horizontal-like map if the first projection
$\pi_1\rest{\Gamma}$ is  injective  outside a locally finite number of curves, and 
$$\overline{\Gamma}\cap\fr(\bb_1\times \bb_2)\subset  (\bb_1\times\fr_v\bb_2 )\cup 
(\fr_h\bb_1\times (\bb_2\cup\fr_v\bb_2)).$$

The crucial fact is that such maps naturally act on families of positive 
closed currents: see proposition \ref{prop_gt} below. 
We need first introduce a few concepts. 

\begin{defi} A horizontal 
positive closed current  in $\bb$ is a  bidegree (1,1)
positive closed current with horizontal  support, that is, 
$$\supp (T)\subset \dd\times \dd_{1-\e} \text{ for some } \e>0.$$
Vertical currents are defined analogously.
\end{defi}

More generally we call {\sl horizontal} an object --form, 
current, or subvariety-- with horizontal
support.
A basic object in this setting is the family of {\sl slice measures} associated to
a  horizontal positive closed current, i.e. the family of measures
$m^z= T\wedge [\set{z}\times\dd]$ supported on the vertical slices $\set{z}\times\dd$. 
These measures do have the same mass (see e.g. \cite{ls, hl}) which is a 
homological-like characteristic of horizontal currents --this can be made precise in
terms of  relative homology. 
We say a horizontal positive closed current is 
{\sl normalized} if its slice mass is 1.\\

There is a useful
canonical plurisubharmonic potential for a horizontal positive closed 
current $T$, given in terms 
of the slice measures. Indeed, let 
\begin{equation}\label{eq_canonical}
u(z,w) = \int_{\set{z}\times\dd} \log\abs{w-s} dm^z(s), 
\end{equation}
$u$ is a  psh function, and $dd^cu=T$. The result is obvious when $T$ 
is the current of integration on an analytic set, hence it holds in the general case 
by approximation of $T$ by horizontal varieties \cite{ds}.   
We call $u$ the {\sl logarithmic} or {\sl canonical potential} of $T$. 

One has the following basic
estimate: if $T$ is normalized, and $\dist(\supp(T), \fr_h\bb)\geq \e$ then 
the $L^\infty$ norm of $u$ in an $\e/2$ neighborhood of $\fr_h\bb$ is bounded 
by $-\log (\e/2)$.

\begin{prop}\label{prop_horiz}
Let $T_1$ and $T_2$ be horizontal positive closed currents of slice mass 1. Then there 
exists a horizontal 1-form $\alpha$ such that $T_1-T_2=d\alpha$. 

Moreover if $T_1$ and $T_2$ are smooth, with respective logarithmic potentials 
$u_1$ and $u_2$, 
and $\supp T_i\subset \dd\times \dd_{1-\e}$, $\alpha$ can be chosen so that
there exists a constant $C$ depending only on $\e$ such that the
following holds 
\begin{equation}\label{eq_alpha}
\norm{\alpha}\leq C(1+\norm{d^c u_1- d^cu_2}),
\end{equation}
where the  norm has to be understood pointwise on coefficients.
\end{prop}

Estimate (\ref{eq_alpha}) is stated for smooth currents, so by the usual 
smoothing process for psh functions, it holds for instance 
when the potentials $u_i$ 
have derivatives in $L^2$. This estimate is important in our setting because 
 bounded psh functions have $L^2$ derivatives (see lemma  \ref{lem_polarization} below).  

\begin{pf} consider the 1-form $d^cu_1-d^cu_2$ in the ``cylinder''
 $\bb\backslash (\dd\times \dd_{1-\e})$; $d^cu_1-d^cu_2$ is closed and the cylinder 
 retracts to any vertical circle $\set{z}\times\fr \dd_{1-\e}$, so the form is exact 
 if and only if its integral on the circle vanishes. 
 Furthermore, 
 $$\int_{\set{z}\times\fr \dd_{1-\e}} d^cu_1-d^cu_2 = \int_{\set{z}\times \dd_{1-\e}} 
 dd^cu_1-dd^cu_2 = 0$$ since the currents have the same slice mass. So
 there exists a function $\theta$ in $\bb\backslash (\dd\times \dd_{1-\e})$ such that 
 $d^cu_1-d^cu_2=d\theta$ there. 
 
 Now let $\chi$ be a cutoff function with horizontal support
 in $\bb$, $\chi=1$ near $\dd\times \dd_{1-
 \frac{\e}{2}}$, and 
 consider the form $$\alpha= \chi (d^cu_1 - d^c u_2) + \theta d\chi,$$
$\alpha$ has horizontal support and one easily checks $d\alpha = T_1-T_2$.\\

To get estimate (\ref{eq_alpha}), note that ${\rm dist}
(\supp(d\chi), \supp(T_i))\geq \e/2$,
so the function $u_1$ and $u_2$ are pluriharmonic and have controlled 
$L^\infty$ norm in the neighborhood of $\supp(d\chi)$. So $d^c u_1- d^cu_2$ is 
a smooth form there, with $\norm{d^cu_1-d^cu_2}_{L^\infty(\supp (d\chi))}\leq C(\e)$.
Since $\theta$ can be constructed by integrating $ d^cu_1-d^cu_2$ along paths, 
we get a control of $\norm{\theta}_{L^\infty(\supp (d\chi))}$ depending
only on $\e$.
\end{pf}

From this one classically deduces the following fact:

\begin{prop} \label{prop_isect}
If $\phi$ is  a vertical closed form and
$T$ is a horizontal positive closed current in $\bb$, 
then $\langle T,\phi\rangle$ only
depends on the 
slice mass of $T$.
\end{prop}

\begin{pf} assume $T_1$ and $T_2$ have the same slice mass, then
  $T_1-T_2=d\alpha$ as
 in the preceding proposition. Since $\supp(\phi)\cap\supp(\alpha)$
 is compact we infer
$$\langle T_1,\phi\rangle- \langle T_2,\phi\rangle= \langle d\alpha, \phi\rangle
=\langle \alpha, d \phi\rangle=0.$$ 
\end{pf}

\noindent{\bf Remark.} 
By smoothing one may replace $\phi$ by any vertical positive closed current $S$, and 
get that $\int T\wedge S$ in computed in relative homology. Indeed the wedge 
product is well defined as a measure in this geometric setting \cite{sib2}.  \\

We now  define the push forward and the pull back Graph
Transform operators for currents associated to a horizontal-like map. This is 
similar to \cite{hl}. Recall that
in dimension 2, for a dominating 
meromorphic map $f$, one can define
 continuous push forward and pull back operators 
$f_*$  and $f^*$ acting on positive closed currents of bidegree (1,1) 
by going through the graph of $f$, i.e. if $\pi_1$, $\pi_2$ denote the projections 
from the graph $\Gamma$, then $f^*T=(\pi_1)_*(\pi_2^*T\wedge [\Gamma])$
and $f_*T=(\pi_2)_*(\pi_1^*T\wedge [\Gamma])$.

\begin{prop} \label{prop_gt}
 Let $f$ be a horizontal-like map  $\bb_1\cv \bb_2$. Then there exists 
an integer $d\geq 1$ such that
\begin{itemize} 
\item[i.] For every normalized horizontal positive closed current $T$ in $\bb_1$, 
$\unsur{d} \mathbf{1}_{\bb_2} f_*T$ is  normalized, horizontal , positive 
and closed in ${\bb_2}$.
\item[ii.] For every normalized vertical positive closed current $T$ in $\bb_2$, 
$\unsur{d} \mathbf{1}_{\bb_1} f^*T$ is  
normalized, vertical, positive and closed in $\bb_1$.
\end{itemize}
\end{prop}

The integer $d$ is by definition the {\sl degree} of $f$. Of course it can be evaluated 
by pushing  a horizontal line in $\bb_1$.

\begin{pf} the push 
forward $f_*$ of a horizontal current $S$ is well defined because 
by the very definition of
horizontal-like maps, for any test 2-form $\phi$ in $\bb_2$, $\supp S\cap
\supp f^*\phi\finc \bb_1$. So if $T$ is smooth and $\varphi$
is a test 1-form, 
 $\langle \mathbf{1}_{\bb_2}
f_*T, d\varphi\rangle= \langle T, d(f^*\varphi)\rangle=0$ since 
$f^*\varphi$ is a vertical form with $L^1_{loc}$ coefficients.

Assume next $T_1$ and $T_2$ are smooth, horizontal, and normalized; 
thus there exists a horizontal
1-form $\alpha$ such that $T_1-T_2=d\alpha$. We can therefore
 write $$ \mathbf{1}_{\bb_2} f_*T_1- \mathbf{1}_{\bb_2}
  f_*T_2 = \mathbf{1}_{\bb_2} f_*d\alpha
= \mathbf{1}_{\bb_2} d f_*\alpha,$$  where $f_*\alpha$ is a horizontal
form with $L^1_{loc}$ 
coefficients. Hence for every vertical closed 2-form $\phi$, 
$$ \langle\mathbf{1}_{\bb_2} f_*T_1- \mathbf{1}_{\bb_2} f_*T_2, \phi\rangle = 
\langle d(f_*\alpha),\phi \rangle =0,$$ which means the two horizontal currents have the 
same slice mass. It is an integer since one may take as $T_1$ a horizontal line.

To get rid of the smoothness hypothesis, just remark that the operator $f_*$ is 
continuous, and that slightly reducing the ambient bidisk if necessary, smoothing 
horizontal currents by convolution still produces horizontal currents.

The case of the pull back $f^*$ is similar.\\

It remains to prove equality of the respective degrees $d(f_*)$ and $d(f^*)$. 
Take a generic
vertical line $L^v$ in $\bb_2$, and a generic horizontal line $L^h$ in $\bb_1$ (in
particular avoiding indeterminacy points); we claim that 
\begin{equation}\label{eq_degrees}
f_*( [L^h] \wedge \mathbf{1}_{\bb_1}f^*[ L^v] ) = \mathbf{1}_{\bb_2}f_*[L^h]\wedge 
[L^v].
\end{equation}
This implies $d(f_*)=d(f^*)$: indeed, pushing forward a measure preserves masses, 
while the mass of  $[L^h] \wedge \mathbf{1}_{\bb_1}f^*[ L^v] $ is $d(f^*)$ and
$\mathbf{1}_{\bb_2}f_*[L^h]\wedge [L^v]$ has mass $d(f_*)$ --observe
that the push forward of 
the measure is well defined because $f\rest{L^h}$ is proper. Now, 
equation  (\ref{eq_degrees}) is an easy consequence of the definitions of $f_*$
and $f^*$.
\end{pf}

The construction in the following section will lead us to consider sequences
of meromorphic horizontal-like maps defined in chains 
of bidisks, and iterate this graph transform 
operator.  Consider a sequence of horizontal-like maps 
$\set{(f_i,\bb_i,\bb_{i+1})}_{i\geq 1}$
of respective degree $d_i$, where the
$\bb_i$ are abstract bidisks, that is, for every $i$,
$\bb_i$ is an open subset of a 2-dimensional 
complex manifold, with a fixed biholomorphism $\varphi_i:\bb\cv\bb_i$ 
--we actually 
assume $\varphi_i$ is defined in a neighborhood of $\bb$. We 
identify $\bb_i$
and $\bb$ through $\varphi_i$, so that we may for instance speak of  
the logarithmic potential of a current in $\bb_i$, etc. \\

We are interested in the following problem. Let 
$\el_i=\unsur{d_i}\mathbf{1}_{\bb_i}
f^* $ be the pull back graph transform on vertical 
currents associated to $f_i$. For each $n$, let  $T_{n+1}$
be a  normalized vertical current in $\bb_{n+1}$, 
  we want to
study the convergence of the sequence of currents $(\el_1\cdots \el_n T_{n+1})$.\\

Our main result is the following. The proof -as well as some more
refined statements-- will be postponed until
section \ref{sec_conv}. 

\begin{thm}\label{thm_cv_main}
Let $\set{(f_i, \bb_i,\bb_{i+1})}_{i\geq 1}$ be a sequence of meromorphic horizontal-like maps, 
such that $f_i^{-1}(\bb_{i+1})\subset \dd_{1-\e}\times \dd\subset \bb_i$ for  
a fixed $\e>0$. Let  $K=\bigcap_{n\geq 1}\big(f_1^{-1}\cdots
f_n^{-1}(\bb_{n+1})\big)$.

Assume for each $n$, $T_{n}$ is a normalized 
positive closed vertical current in $\bb_{n}$, with 
 canonical potential uniformly 
 bounded with respect to $n$. Assume also either $K$ is polynomially
 convex or has zero Lebesgue measure.
Then the sequence of currents 
$(\el_1\cdots \el_n T_{n+1})$ converges to a current $\tau$ in $\bb_1$, independent 
of the sequence $(T_n)$. 
\end{thm}

Notice that in case the maps $f_i$ are holomorphic $\bb_i\cv \cd$, the
polynomial convexity assumption is fulfilled. 
With more restrictive hypotheses on the maps $f_i$, we are able 
to say more about the current $\tau$. Indeed assume the $f_i$
are holomorphic outside the origin, with a possible indeterminacy point
at 0. Assume the orbits approach the origin 
in a somewhat tame fashion: 
\begin{equation}\label{eq_attract_both}
\forall i,~\forall p=(z,w)\in f^{-1}_i(\bb_{i+1}),~
c\abs{w}^{\ell_i} \leq \dist (f_i(p),(w=0)) \leq \abs{w}^{\ell_i}.
\end{equation} 
We also suppose for convenience 
that the sequence of maps is an infinite  word in a finite 
number of ``letters" $f_i$. These assumptions are satisfied for the
mappings considered in the next section.

Then 
\begin{itemize}
\item[-] if $\displaystyle{\frac{\ell_1\cdots\ell_n}
{d_1\cdots d_n}}$ tends to infinity, $\tau$ has pluripolar support
(proposition \ref{prop_polar}), 
\item[-] if the series $\displaystyle{\sum_{n\geq 1} \frac{1+n\ell_1\cdots\ell_n}
{d_1\cdots d_n}}$ converges, $\tau$ has continuous 
 potential out of the
line $(w=0)$ (proposition \ref{prop_apriori}).
\end{itemize}


\section{The general picture}\label{sec_speed}

In this section we present the class of polynomial maps in $\cd$ that 
motivates our study. A polynomial map in $\cd$ can be seen as a meromorphic map
on a  compactification of $\cd$, which is usually chosen to ensure ``algebraic 
stability" if possible. For simplicity we restrict to the compactification
$\pd$. On $\pd$ it is 
known how to construct a dynamically interesting invariant current
 under the algebraic stability hypothesis \cite{sib}.\\

A good starting point for the finer study of the dynamics is to analyze the behaviour 
of the invariant current $T$ near an invariant set where the dynamics 
is known. This has been considered by J.E. Fornaess and the third author 
\cite{fs} for a class of holomorphic endomorphisms of $\pd$.
 Bedford and Jonsson \cite{bj} 
considered the case of a polynomial map $f$ in $\cc^k$ with holomorphic continuation to
$\mathbb{P}^k$. 
They explore the interplay between the 
structure of the Green current $T$
in the neighborhood of infinity in $\cc^k$ 
and the dynamics on the hyperplane at infinity $H_\infty$. 
It appears that $T^{k-1}$ admits a local laminar decomposition, consisting of
local stable manifolds of the maximal entropy measure of $f\rest{H_\infty}$.\\

Here we consider the opposite situation, in dimension 2 only. 
There is no interesting dynamics on the 
line at infinity, whereas some complexity 
arises from the indeterminacy points. The Green current admits a decomposition into
pieces that reflect the attracting behaviour of indeterminacy points, but since the image 
of each of these is the whole line at infinity, they cannot be treated separately. 
This feature is overcome by using partitions and symbolic dynamics.
The study of the mixed case where there are both non trivial dynamics at infinity
and attracting indeterminacy points remains to be done as well as the case where the 
divisor at infinity has several irreducible components. 

\subsection{The speed spectrum}\label{subs_speed}
Let us  be more specific. Let $f$ be a 
polynomial map of degree $d$ in $\cd$ such that infinity
is $f$-attracting, that is, for $\norm{p}$ large enough $\norm{f(p)}\geq c\norm{p}$
for some $c>1$.  We assume (the rational 
extension of)
$f$ has $m\geq 1$ indeterminacy points $I_1,\ldots, I_m$ on the line at infinity 
$L_\infty\subset\pd$ and $f(L_\infty\backslash I(f))=X\notin I(f)$ is
a single
 fixed point. 
Hence $f$ is algebraically stable in $\pd$, i.e. for every $n\geq 1$,
$\mathrm{deg}(f^n)=(\mathrm{deg}f)^n$.
 Since the image of each indeterminacy point 
must be a non trivial variety, not intersecting $\cd$ because $f$ is proper, one has
$f(I_k)=L_\infty$ for each $k$. Assume $X=[1:0:0]$; this means that
we performed a linear change of coordinates so that $f=(f_1,f_2)$ with 
 $d=d_1=\mathrm{deg}(f_1)>\mathrm{deg} (f_2)=d_2$.
We fix a neighborhood $V(X)$ of the point $X$,
of the form $V(X)=\set{\abs{w}\leq \e\abs{z}, ~\abs{z}\geq 1/\e}$,  
such that $I(f)\cap V(X)=\emptyset$ and $f(V(X))\finc V(X)$. \\

We also need a precise control on the rate of attraction of points to the line 
$L_\infty$. This is well understood near $X$ where the points go to infinity at maximum
speed $\norm{p}^{d^n}$. 
When $f$ is {\sl regular}, that is $\set{f_1^+=f_2^+=0}$ is reduced to zero, where 
$f^+$ denotes the homogeneous term of highest degree (in other terms 
$f_2^+$ does not vanish on $I$), then out of the basin of $X$ one has 
$C_1\norm{p}^{d_2}\leq \norm{f(p)}\leq C_2\norm{p}^{d_2}$. The more general class 
of {\sl semi-regular} 
maps is also studied in \cite{dis2}. Besides these cases
there is 
no such estimate near the indeterminacy points
so we need to
formulate a hypothesis. More precisely we suppose that  for every $i$,
there exist constants $C_1, C_2$, 
  ${\ell}_i\geq 1$, and neighborhoods $V(I_i)$ such that 
\begin{equation}\label{eq_speed}
\forall p\in V(I_i), \text{ s.t. } f(p)\notin V(X),~C_1\norm{p}^{{\ell}_i}
\leq \norm{f(p)} \leq C_2\norm{p}^{{\ell}_i}.
\end{equation}
Notice that if $\ell_i=1$, since $L_\infty$ is attracting, one has
$C_1>1$. For  $\ell_i>1$, performing a linear change of coordinates if
necessary,
one may always assume $C_1\geq 1$.
We will compute these exponents on some examples in \S
\ref{subs_examples}, and in a more systematic way in proposition 
\ref{prop_newton}.
In some cases it could be more convenient to consider upper and lower
 exponents $\overline{\ell}_i$ and  $\underline{\ell}_i$ in each $V(I_i)$. \\

We define the {\sl lower} and {\sl upper escape rates} (or
{\sl divergence speed}) of the point 
$p$ respectively as
$$\log \underline{\ell}(p)=
\liminf_{n\cv\infty} \unsur{n}\log^+ \log^+{\norm{f^n(p)}} \text{ and }
\log \overline{\ell}(p)=
 \limsup_{n\cv\infty}  \unsur{n}\log^+ \log^+{\norm{f^n(p)}}.$$ Of course when 
the two coincide we simply speak of  {\sl escape rate} $\ell(p)$.
In the case  of semi-regular maps, 
the function $\ell$ takes  only finitely many values, see
\cite{dis2}.
 The new phenomenon here
  is that when for some $i,j$,  $\ell_i\neq\ell_j$,  there is a whole interval 
of possible speeds of divergence to infinity, which is in sharp
contrast with all previously known examples. We will exhibit 
mappings for which the interval $[\min \ell_i, \max \ell_i]$ is non trivial
in \S\ref{subs_examples} below.

\begin{thm}\label{thm_speed}
Assume $f$ is a polynomial map in $\cd$ of degree $d$, satisfying (\ref{eq_speed}),
such that $L_\infty$ is attracting 
and  $f(L_\infty\backslash I(f))=X\in L_\infty$.
Then, with notations as before,
 the spectrum of possible escape rates of points in $\cd$ to infinity
is exactly $[\min \ell_i, \max \ell_i]\cup\set{d}$, with $\max \ell_i\leq d-1$.

Moreover a point has escape rate 
$d$ if and only if it is in the basin of $X$.
\end{thm}

The proof will follow from a sequence  of observations.
The basic geometric feature is the following.

\begin{prop}\label{prop_hlinfty}
There exists a family of bidisks $\bb_i\ni I_i$, $i=1, \ldots, m$, 
so that for every pair $(i,j)$, $(f,\bb_i, \bb_j)$ is a horizontal-like map, 
of degree $d_i$ depending only on $i$. Moreover $\sum_i d_i=d$.
\end{prop}

\begin{pf} write the coordinates in $\cd$ as $(z,w)=(\frac{Z}{T}, \frac{W}{T})$, where
$[Z:W:T]$ are the homogeneous coordinates in $\pd$. Since $X=[1:0:0]$, near 
$I_i$ one may 
use the  coordinates $$(u,v)=(\frac{z}{w}, \unsur{w})=(\frac{Z}{W},
\frac{T}{W}), $$  
in which the line at infinity becomes  $(v=0)$. Our bidisk neighborhoods have the 
form $$\bb_i=\bb_{i,r,r'}=\set{\abs{u-u_i}\leq r,~ \abs{v}\leq r'}\text{ , 
where }I_i=(u_i,0), $$
where  $r$ is chosen
so small that the bidisks are mutually disjoint and disjoint from 
$V(X)$. The radius $r$ being fixed, since  $f(L_\infty\backslash I(f))=X$, 
for $r'$ small enough, $f(\fr_v\bb_i)\subset V(X)$. 
That $f(\bb_i)\cap\bb_j\neq \emptyset$ for every $(i,j)$ follows from
the fact that $f(I_i)=L_\infty$; moreover, since $L_\infty$ 
is attracting, $f(\bb_i)$ is horizontally supported in $\bb_j$. Thus the hypotheses of 
definition \ref{def_hl} are satisfied.\\

On the other hand, if $L_i$ is a perturbation of $L_\infty$ of the form
 $(v=c), ~c\neq 0$ in $\bb_i$, 
$f(L_i)$ is a curve in $\cd$ with boundary in $V(X)$, close to $L_\infty$,
 so the linear projection 
$\pi_0:\cd\ni (z,w)\mapsto[z:w]\in L_\infty$ makes it a branched cover of some degree 
$d_i$ over $L_\infty\backslash V(X)$. The integer $d_i$ is thus the degree of 
$(f,\bb_i,\bb_j)$ for every $j$. Of course viewing the line $(v=c)$ as a global 
line in $\pd$, $f(L)$ has algebraic degree $d$, which coincides with the degree
of the branched cover $\pi_0\rest{f(L)}$. Over $L_\infty\backslash V(X)$, the 
contributions come from the boxes  $\bb_i$, so $\sum d_i=d$.
\end{pf}

As a consequence of the proof one remarks that 
if $L'$ is a generic line, 
$d_i$ is exactly the intersection
multiplicity of $f^{-1}(L')$ and $L_\infty$ at $I_i$, since a generic perturbation $L$ of
$L_\infty$ intersects  $f^{-1}(L')$ in $d_i$ points in $\bb_i$. 

\begin{prop}\label{prop_formula} 
The topological degree of $f$ satisfies
$\displaystyle{d_t=\sum_{i=1}^m \ell_i d_i}$. In particular $d_t\geq d$.
\end{prop}

\begin{pf} let $L$ be a generic line in $\pd$, avoiding $V(X)$. 
On $L$ and $f^{-1}(L)$ one
may consider the coordinates $(u,v)$ as in the previous proof. The psh potential of
$(v=0)=L_\infty$ is $\log\abs{v}$, and one obviously has
$$\int dd^c\log\abs{v}\wedge [L] = \int [L_\infty]\wedge [L] =1.$$ We will compute 
$\int f^*(dd^c\log\abs{v}\wedge [L])$ in two ways. 

First, note that the integral is supported in a neighborhood of
$L_\infty$
and replace $\log\abs{v}$ by a decreasing approximation
$\varphi_\e=\max(\log\abs{v}, \log\e)$. 
The pull back of the probability measure
$dd^c\varphi_\e\wedge [L]$  by $f$ has mass $d_t$ by definition of the topological 
degree. 

On the other hand the integral 
$$\int f^*(dd^c\varphi_\e\wedge[L])= \int dd^c(\varphi_\e\rond f) \wedge [f^{-1}(L)]
=\int_{f^{-1}(L)} dd^c(\varphi_\e\rond f)$$ may be computed directly.
Indeed $\varphi_\e\rond f(p) = \max(\log |{\tilde f_2(p)} |,\log{\e})$, where $\tilde f_2$ 
denotes
the second coordinate, with respect to the coordinate system
$(u,v)$. If $p=(u,v)\in f^{-1}(L)\cap\bb_i$, 
$$ \ell_i\log\abs{v} +O(1)\leq \log |{\tilde f_2(p)} |\leq \ell_i\log\abs{v}
+O(1) $$ by  definition of 
$\ell_i$ --remember both $L$ and $f^{-1}(L)$ are far from $X=[1:0:0]$ so 
the coordinate $v$ determines $\norm{p}$ in $\cd$. Hence $\varphi_\e\rond f$ decreases
to $\log{ |{\tilde f_2} |}=\ell_i\log\abs{v}+O(1)$ in $f^{-1}(L)\cap\bb_i$. 
Since for small $\e$, $\varphi_\e\rond f\rest{f^{-1}(L)}$ 
is harmonic outside a small neighborhood of $I(f)=\set{I_1,\ldots, I_m}$, we infer
$$\int_{f^{-1}(L)} dd^c(\varphi_\e\rond f) = \sum_i \int_{f^{-1}(L)\cap\bb_i} 
dd^c(\varphi_\e\rond f) \underset{\e\cv 0}{\longrightarrow} 
\sum_i \int_{f^{-1}(L)\cap \bb_i}dd^c(\ell_i\log\abs{v}+O(1))$$
where the $dd^c$ vanishes outside 
 any arbitrarily small neighborhood of $I_i$. This gives
$$\int f^*(dd^c\varphi_\e\wedge [L]) \longrightarrow
 \sum_i \ell_i \int[L_\infty]\wedge [f^{-1}(L)]
=\sum_i\ell_i d_i$$ since the intersection multiplicity of $f^{-1}(L)$ and $L_\infty$ at
$I_i$ is $d_i$.
\end{pf}

Let $\mathcal{K}$ be the closure in $\pd$ of the 
complement of the basin of attraction of $X$. We put $N=\bb_1\cup\cdots \cup \bb_m$,
so that in a neighborhood of 
 $L_\infty$ one has $\mathcal{K}\subset N $.
 Since the line $L_\infty$ is attracting, this allows 
to introduce symbolic dynamics in $\mathcal{K}\cap N$. More precisely for 
any sequence $\alpha=(\alpha(j))_{j\in\nn}\in
\set{1,\ldots ,m}^\nn$, let 
$$\mathcal{K}_\alpha= \set{p\in \mathcal{K}\cap N, ~f^j(p)\in \bb_{\alpha(j)}}$$
be the set of points with itinerary $\alpha$. This set is in a sense the local
stable set of the sequence $\alpha$. 
By proposition \ref{prop_hlinfty} and the discussion in section \ref{sec_hl}, 
for every $\alpha$, $\mathcal{K}_\alpha\subset\bb_{\alpha(0)}$ 
is a non empty vertical closed set, carrying a closed positive current. We will 
discuss convergence properties of the sequence of currents in the next paragraph.
Before this, we compute the escape rate of points on $\mathcal{K}_\alpha$, through
the following lemma.

\begin{lem}\label{lem_speed}
Let $\alpha\in\set{1,\ldots ,m}^\nn$. 
\begin{enumerate}[\em i.]
\item For  fixed $p$ in $\mathcal{K}_\alpha$,
$$\unsur{n}\log\log \norm{f^n(p)}= \unsur{n}\log( \ell_{\alpha(0)}
\cdots\ell_{\alpha(n-1)}) 
+o(1)\text{ as } n\cv\infty.$$
\item If moreover the series 
${\sum_{n\geq 1} \unsur{\ell_{\alpha(0)}\cdots\ell_{\alpha(n)}}}$ converges, the 
sequence of functions $$\unsur{\ell_{\alpha(0)}\cdots\ell_{\alpha(n-1)} }
\log\norm{f^n(p)}$$ converges uniformly to a function $G_\alpha$ on 
$\mathcal{K}_\alpha$.
\end{enumerate}
\end{lem}

Note that the hypothesis in the second item is obviously satisfied when all $
\ell_i\geq 2$. We do  not know whether  $\mathcal{K}_\alpha$ has laminar
structure in general, nevertheless one easily sees  $G_\alpha$ is harmonic 
on any holomorphic disk $\Delta\subset \mathcal{K}_\alpha$, if any.

\begin{pf} if $p\in\mathcal{K}_\alpha$,
by (\ref{eq_speed}), for every $0\leq j\leq k-1$ one has
--recall both $C_1, C_2\geq 1$--
$$\log C_1\leq \log\norm{f^{j+1}(p)}-\ell_{\alpha(j)}\log\norm{f^j(p)}\leq \log C_2.$$
Summing this formula from $j=0$ to $k-1$ yields 
$$
 \log\|{f^{k}(p)}\|
-\ell_{\alpha(0)}\cdots\ell_{\alpha(k-1)}\log\norm{p} \leq 
\left(\sum_{j=0}^{k-1}  {\ell_{\alpha(j+1)}\cdots\ell_{\alpha(k-1)} } \right)
 \log C_2   ,$$
with a similar left hand side inequality involving $\log C_1$, 
so 
$$\abs{\unsur{k}\log\log \|{f^k(p)}\|- \unsur{k}\log(
\ell_{\alpha(0)}\cdots\ell_{\alpha(k-1)})} \leq \unsur{k}
\log\left( \log\norm{p} + C\sum_{j=0}^{k-1}
  \unsur{\ell_{\alpha(0)}\cdots\ell_{\alpha(j)}} \right) .$$
Since the sum in the right hand side is $O(k)$,
we get item i.

To prove the second point, just remark that 
$$\abs{\unsur{\ell_{\alpha(0)}\cdots\ell_{\alpha(j)} }
\log\norm{f^{j+1}(p)}- \unsur{\ell_{\alpha(0)}\cdots\ell_{\alpha(j-1)} }
\log\norm{f^j(p)}} \leq \frac{
\max(\abs{ \log C_1}, \abs{\log C_2})}{\ell_{\alpha(0)}\cdots\ell_{\alpha(j)} }$$
on $\mathcal{K}_\alpha$, and apply Cauchy's criterion.
\end{pf}

\noindent{\bf Proof of theorem \ref{thm_speed}:}
 from the previous lemma we  
 know all possible speeds of divergence to infinity. It is clear that for 
every $\ell\in[\min \ell_i, \max \ell_i]$ there exists a sequence $\alpha$ such that
$\ell=\lim (\ell_{\alpha(0)}\cdots\ell_{\alpha(k-1)})^{\unsur{k}}$. On the other hand 
$\max \ell_i \leq d-1$: indeed assume $f=(f_1,f_2)$, with $d=\mathrm{deg}(f_1)>
\mathrm{deg}(f_2)$ as before, then if $\norm{p}$ is large and $f(p)$ is 
far from $X$, then $\abs{f_2(p)}\geq c \abs{f_1(p)}$. This means 
that $\norm{f(p)}\leq C\norm{p}^{\mathrm{deg}(f_2)}$. This concludes the proof. 
 \hfill$\square$

\subsection{Decomposition of the Green current.}\label{subs_decompos}

In this section we show how the structure of the Green current 
 is determined by the above picture. Near infinity, to each itinerary sequence 
 $\alpha$ corresponds a unique current $T_\alpha$, which is part of the Green
 current $T$. This Cantor structure persists after finitely many backward iterations,
 modulo some  identifications between itineraries. 
 
Consider the shift space $\Sigma=\set{1,\ldots,m}^\nn$, and let $\nu$ be the 
self similar
probability measure on $\Sigma$ defined by 
$$\nu\left(\set{\set{i}\times\set{1,\ldots,m}^{\nn^*}}\right)=\frac{d_i}{d},$$
that is, the measure of the cylinder of sequences starting with
$\alpha(0), \ldots, \alpha(k)$ is $\unsur{d^{k+1}}d_{\alpha(0)}\cdots
d_{\alpha(k)}$. Under the action of the usual one sided shift $\sigma$ 
on $\Sigma$, $\nu$ is invariant, and  mixing. 

We keep the hypotheses of the preceding section, namely we
assume the polynomial map $f$ has $m\geq 2$ indeterminacy points
at infinity, and the
numbers $d_i$  have the same meaning. Recall that for
every pair $(i,j)$, $(f, \bb_i, \bb_j)$ defines a horizontal-like map,
hence acts on vertical closed positive currents with degree $d_i$. 
We denote by $\el_i= \unsur{d_i}\mathbf{1}_{\bb_i}f^*$ 
the corresponding pull back operator, 
which can be viewed as being independent of $j$. Notice that
hypothesis (\ref{eq_speed})
is unnecessary in the next theorem.

\begin{thm}\label{thm_decompos}
Under the preceding hypotheses, the following is true:
\begin{enumerate}[\em i.]
\item There exists an at most countable set
  $\mathcal{E}\subset\Sigma$,
such that 
for  every sequence $\alpha\in\Sigma\setminus\mathcal{E}$ 
there exists a
  current $T_\alpha$ supported in $\mathcal{K}_\alpha$  uniquely
  defined as follows: for every vertical positive closed current
  $S$ in $\bb_1\cup\cdots\cup\bb_m$ with bounded potential
$$\el_{\alpha(0)}\cdots\el_{\alpha(k)}S \cvf T_\alpha.$$
The current $T_\alpha$ is a normalized vertical positive current in 
$\bb_{\alpha(0)}$ and 
satisfies the compatibility relation $\el_iT_\alpha=
T_{i\alpha}$, where $i\alpha$ is the concatenated sequence. 
\item The Green current $T$ admits the decomposition 
$\displaystyle{T=\int_\Sigma T_\alpha d\nu(\alpha)}$ in a neighborhood of infinity. 
\end{enumerate}
\end{thm}
 
\begin{cor}\label{cor_decompos} Assume moreover $f$ satisfies (\ref{eq_speed}).
Then:
\begin{enumerate}[\em i.] 
\item For $\norm{T}$ almost every point $p$ in the neighborhood
of infinity, the escape rate $\ell(p)$ is 
$$\ell(p)= \prod_{i=1}^m \ell_i^\frac{d_i}{d}.$$ 
\item If $\displaystyle{ \prod_{i=1}^m\left(\frac{\ell_i}{d_i}\right)^{\frac{d_i}{d}}>1}$ then
for $\nu$-a.e. $\alpha$, $\mathcal{K}_\alpha$ is pluripolar.
\item If $\displaystyle{ \prod_{i=1}^m\left(\frac{\ell_i}{d_i}\right)^{\frac{d_i}{d}}<1}$ then
for $\nu$-a.e. $\alpha$, $T_\alpha$ has continuous potential outside $L_\infty$.
\end{enumerate}
\end{cor}

Item i. in the corollary still holds for almost every $p\in \supp(T)$ escaping to
infinity since $\ell(p)$ does not vary under finitely many iterates. 
Notice that concavity of the logarithm function together with proposition
\ref{prop_formula} imply that the
generic speed $\ell(p)$ is not larger than $d_t/d$, with strict
inequality when at least two $\ell_i$'s are different.

Item ii. and iii. hold under
slightly weaker hypotheses: suppose for every $i$,
there exist constants $C_1, C_2$, 
  $\underline{\ell}_i\geq 1$ and   $\overline{\ell}_i\geq 1$ such that 
\begin{equation}\label{eq_speedmodified}
\forall p\in\bb_i, \text{ s.t. } f(p)\notin V(X),~C_1\norm{p}^{\underline{\ell}_i}
\leq \norm{f(p)} \leq C_2\norm{p}^{\overline{\ell}_i}.
\end{equation}
Then in ii. (resp. iii.) $\ell_i$ may be replaced by $\underline{\ell}_i$ 
(resp. $\overline{\ell}_i$). \\

\noindent{\bf Proof of the theorem:} if $S$ is a smooth
vertical closed positive  current in
$\bb_1\cup\cdots\cup\bb_m$, normalized in the sense that its intersection
with a global horizontal line has mass  1, and with bounded potential
$$\unsur{d} f^*S= \sum_{i=1}^m \frac{d_i}{d}\el_i S,$$ so 
$$\unsur{d^k} (f^k)^*S = \sum_{\set{\alpha(0),\ldots,\alpha(k-1)}\in
  \set{1,\ldots,m}^k} \frac{d_{\alpha(0)}\cdots 
d_{\alpha(k-1)}}{d^k} \el_{\alpha(0)}\cdots\el_{\alpha(k-1)} S\cvf T. $$  
In the neighborhood of infinity, $\mathcal{K}$ is the disjoint union
  of the sets $\mathcal{K}_\alpha$, so for $\alpha$ not in  an
  at most countable set $\mathcal{E}$, $\mathcal{K}_\alpha$ has zero Lebesgue
  measure. 
For every such sequence $\alpha$,  
$ \el_{\alpha(0)}\cdots\el_{\alpha(k-1)} S$ 
converges by theorem \ref{thm_cv_main}
and  the
decomposition formula $T=\int T_\alpha d\nu(\alpha)$
follows from  the dominated convergence theorem. \hfill $\square$\\

\noindent{\bf Proof of the corollary:}
the itinerary map
$\mathcal{K}\cap N\cv\Sigma$ semi conjugates $f$ to the shift $\sigma$
on $\Sigma$,
i.e. $f(\mathcal{K}_\alpha)\subset\mathcal{K}_{\sigma(\alpha)}$. On
the other hand $\sigma$ is $\nu$-ergodic, so by the Birkhoff ergodic
theorem, for $\nu$-a.e. $\alpha$ one has
$$\left({\ell_{\alpha(0)} 
    \cdots\ell_{\alpha(k-1)}}\right)^{\unsur{k}}= \exp\left(
\unsur{k}\sum_{i=0}^{k-1}\log\ell_{\sigma^i\alpha(0)} \right)\underset{k\cv\infty}{\longrightarrow}
\prod_{i=1}^m\ell_i^{\frac{d_i}{d}}= \exp\left(
    \int_\Sigma \log
    \ell_{\alpha(0)}d\nu(\alpha) \right). $$ 
Thus the first item of the corollary follows from  lemma \ref{lem_speed}
and the decomposition of the Green current.
    
In the same vein, using the Birkhoff ergodic theorem for the sequence 
$$\left(\frac{\ell_{\alpha(0)} 
    \cdots\ell_{\alpha(k-1)} } {d_{\alpha(0)}\cdots d_{\alpha(k-1)}}
    \right)^{\unsur{k}}$$ 
    shows that:
 \begin{itemize}
 \item[-] either $\prod_{i=1}^m\left(\frac{\ell_i}{d_i}
\right)^{\frac{d_i}{d}}>1$ then for generic $\alpha$, $
\displaystyle{\frac{\ell_{\alpha(0)}\cdots
\ell_{\alpha(n)} }{d_{\alpha(0)}\cdots d_{\alpha(n)} }\cv\infty}$, and proposition
\ref{prop_polar} gives ii.;
\item[-] or $\prod_{i=1}^m\left(\frac{\ell_i}{d_i}
\right)^{\frac{d_i}{d}}<1$ then for generic $\alpha$
the series $\displaystyle{\sum \frac{1+n\ell_{\alpha(0)}\cdots
\ell_{\alpha(n)} }{d_{\alpha(0)}\cdots d_{\alpha(n)} }}$ converges, 
and iii. follows from proposition \ref{prop_apriori}.
\end{itemize}
\hfill $\square$\\

There is a alternate decomposition of the Green current, related to the spectrum 
of escape rates. Consider for example
$$\log \underline{\ell}(p)=
\liminf_{n\cv\infty} \unsur{n}\log \log{\norm{f^n(p)}}= \liminf_{n\cv\infty}
\unsur{n}\log (\ell_{\alpha(0)}\cdots\ell_{\alpha(n-1)}) \text{ where }p\in
\mathcal{K}_\alpha.$$
This defines a  measurable function 
$\alpha \mapsto  \liminf\unsur{n}\log (\ell_{\alpha(0)}\cdots
\ell_{\alpha(n-1)}) $  on $\Sigma$, whose level sets  form a measurable 
partition --most level sets have zero measure since there exists a generic speed.
The measure
 $\nu$ admits conditional 
probability measures on almost every
 level set $\set{\underline{\ell}=c}$, and we deduce that in a neighborhood of
infinity, almost every level set $\set{\underline{\ell}=c}$ 
carries a natural positive closed current defined by 
$$\int T_\alpha d\nu(\alpha\vert\underline{\ell}=c),$$  
where $\nu(\cdot\vert\underline{\ell}=c)$ denotes conditional
measure.\\

We now explain how  the decomposition 
of theorem \ref{thm_decompos} can be extended to the basin of attraction 
$\mathcal{K}'$ of $L_\infty$ in $\mathcal{K}$ --the complement in
$\mathcal{K}$ of the set of points with bounded orbits. 
Let $N=\bb_1\cup\cdots\cup\bb_m=  
N_0$ and $N_k=f^{-k}(N)$, so that
$\bigcap_{j\geq 0}
\bigcup_{k\geq j} N_k=\mathcal{K}'$.  The point is that even if the currents $T_\alpha$
themselves need not extend to $N_k$, some averaged $T_{{\alpha}^{(k)}}$ do. 
As $k\cv\infty$, more and more 
information  is lost on the decomposition. Recall most preimages of points
converge to the maximal entropy measure as $k\cv\infty$, so the decomposition 
probably has to do with the structure of the measure itself.

Given $\alpha\in\Sigma $, we let ${\alpha}^{(k)}$ be the sets of its $k$-neighbors,
i.e.  ${\alpha}^{(k)}=\sigma^{-k}\sigma^k(\alpha)$, and $\Sigma^{(k)}$ denotes 
the quotient of $\Sigma$ by the $k$-neighborhood relation.  $\Sigma^{(k)}$ is 
isomorphic to $\Sigma$. 

\begin{prop}\label{prop_extend}
Under the hypotheses of theorem \ref{thm_decompos}, for every $k\geq 0$ one has:
\begin{enumerate}[\em i.]
\item The normalized
current $\displaystyle{T_{{\alpha}^{(k)}}=\sum_{\beta\in{\alpha}^{(k)}} 
\frac{d_{\beta(0)}\cdots d_{\beta(k-1)}}{d^k} T_\beta}$ admits an extension to $N_k$.
\item There exists a probability measure $\nu^{(k)}$ on
$\Sigma^{(k)}$ such that the Green current  
decomposes as $\displaystyle{ T=\int_{\Sigma^{(k)} }
  T_{{\alpha}^{(k)}}d\nu^{(k)}}$ on 
$N_k$.
\end{enumerate}
\end{prop}

\begin{pf} the proof is easy. We define the extension by $T_{\alpha^{(k)}}= \unsur{d^k}
(f^{k})^* T_{\sigma^k(\alpha)}$. The current $T_{\alpha^{(k)}}$ is well defined in 
$f^{-k}(N)$ and coincides with the given expression in $N$.

Moreover using the invariance of $\nu$ under $\sigma$, one may write 
$T=\int_\Sigma T_{\sigma^k(\alpha)}d\nu(\alpha)$ in $N$. 
Pulling back this expression
yields
$$T\rest{f^{-k}N} = \unsur{d^k} (f^k)^*(T\rest{N})=\int_\Sigma \unsur{d^k} (f^k)^*
(T_{\sigma^k(\alpha)}) d\nu(\alpha) =\int_\Sigma T_{\alpha^{(k)} }d\nu(\alpha). $$
Since $T_{\alpha^{(k)}}$  only depends on $\alpha^{(k)}\in\Sigma^{(k)}$ we 
get the desired measure on the quotient space $\Sigma^{(k)}$.
\end{pf}

\subsection{Examples.}\label{subs_examples}

The following basic class of examples was the first motivation for this
study.
\begin{example} Let $f:\cd\cv\cd$ be the polynomial map defined by
$$f(z,w)=( z^c(w-z)^d, w^a+z^b).$$
We moreover assume $c+d>b\geq a\geq2$.
Then $f$ has degree
$c+d$, the indeterminacy set is $I=\set{I_1=[0:1:0],I_2=[1:1:0]}$
in homogeneous coordinates $[z:w:t]$, whereas $f(L_\infty\backslash I)=X=[1:0:0]$. 
When $a=b$, $f$ is regular in the sense of \cite{dis2}, so the first new case
is $b>a$.
We will prove condition (\ref{eq_speed}) holds and
 $d_1=c$, $d_2=d$, $\ell_1=a$, and 
$\ell_2=b$. Since $a\geq 2$ $f$ is proper and the line at infinity is attracting.
This example shows both cases ii. and iii. of corollary
\ref{cor_decompos}  may occur. 

Computing the $d_i$ is quite simple. Take a generic line, say $z=\alpha w+\beta$,
and consider its preimage by $f$, which is the curve of equation
$z^c(w-z)^d=\alpha(w^a+z^b)+ \beta$. This curve intersects $L_\infty$ with 
multiplicity $c$ at $I_1$ and $d$ at $I_2$.  

Let us estimate $\ell_1$. For ease of reading we use the following notation: 
$x\lesssim y$ means $x\leq Cy$ for some constant $C$, and $x\ll y$ means 
$x=o(y)$ at infinity. If $p=(z,w)$ is near $I_1$ and $f(p)$ is far from 
$X$ then $\abs{z}\leq\abs{w}/2$, say, and $\abs{w^a+z^b}\gtrsim\abs{z^c(w-z)^d}$.
From the first relation one deduces 
$$\abs{w-z}^d\gtrsim\abs{w}^d\gtrsim\abs{z}^d,$$
and  from the second one, 
$$\abs{w}^a+\abs{z}^b\geq \abs{w^a+z^b}\gtrsim\abs{z^c(w-
z)^d}\gtrsim\abs{z}^{c+d}.$$
Since $c+d>b$, this implies $\abs{z}^{c+d}\lesssim \abs{w}^a$, i.e. $\abs{z}\lesssim
\abs{w}^{\frac{a}{c+d}}$. Thus 
$$\abs{z}^b\lesssim\abs{w}^{\frac{ab}{c+d}}\ll\abs{w}^a,$$ for $\frac{b}{c+d}<1$.  We 
conclude that 
$$ \abs{w}^a\lesssim \abs{w}^a-\abs{z}^b \leq \abs{w^a+z^b} \leq 
 \abs{w}^a+\abs{z}^b\lesssim \abs{w}^a,$$ that is, $\ell_1=a$.
 
 The computation for $\ell_2$ is of the same kind: near $I_2$ one has 
 $\abs{z}\lesssim\abs{w}\lesssim\abs{z}$, so $\abs{w}^a\lesssim\abs{z}^a
 \ll\abs{z}^b$. We infer 
  $$\abs{z}^b\lesssim \abs{z}^b-\abs{w}^a\leq \abs{w^a+z^b} \leq 
 \abs{w}^a+\abs{z}^b\lesssim \abs{z}^b,$$ so $\ell_2=b$.\hfill $\square$
 \end{example}

\begin{example} The following example shows a larger number of indeterminacy 
points is possible. Let 
$$f(z,w)=(z^{d_1}(z-w)^{d_2}(z+w)^{d_3}, z^a+z^b(z-w)^c+ (z+w)^d), $$
where $c<d<a<b+c$, $b<a$, and large $d_2$. 
We leave as an exercise to the reader to prove 
$I_1=[0:1:0]$, $I_2=[1:1:0]$, $I_3=[1:-1:0]$, and 
$\ell_1=d$, $\ell_2=a$, $\ell_3=b+c$. \hfill $\square$
\end{example}

More generally, it is possible  to compute 
the exponents using  {\sl Newton polygons}. This is similar in
spirit to  \cite[\S 3]{dis2}. Here we compute the exponent at 
$I_1=[0:1:0]$, nevertheless a linear change of coordinates allows to 
treat other points on $L_\infty$ as well.

Let $f=(f_1,f_2)$ be a polynomial map with 
$\deg f_1>\deg f_2$. Define the Newton polygon
$\Pi_i$ to be the set of $(r,s)\in\nn^2$
such that the coefficient of $z^rw^s$ in $f_i$ is non zero.

\begin{prop} \label{prop_newton}
Assume there exist  positive integers 
$r_0$, $r_1$, $s_0$, $s_1$, with $0<s_0<r_0=r_1+s_1$, such that
$(r_0,0)\in \Pi_1$, $(0,s_0)\in\Pi_2$, $(r_1,s_1)\in\Pi_1$ and 
$\Pi_1\cup\Pi_2$ is contained in the closed $4$-gon $\mathcal{P}$ 
with vertices $(0,0)$, 
$(r_0,0)$, $(r_1,s_1)$ and $(0,s_0)$. Assume also 
that $\Pi_2\setminus (0,s_0)$
does not meet the segment $\mathcal{S}$ joining $(0,s_0)$ and $(r_1,s_1)$.
Then $I_1=[0:1:0]$ is a point 
of indeterminacy of $f$ with  escape rate $s_0$. 
\end{prop}

\begin{pf} it is clear that $\deg f_1=r_0$ and if $f_1^+$ is the homogeneous 
part of maximal degree of $f_1$ then $f_1^+(0,1)=0$. It follows that $I_1$
is a point of indeterminacy.

Let $p_i=(z_i,w_i)$ be a sequence of points converging to $I$ such that 
$f(p_i)\not\in V(X)$. We need only to check that $|f_2(p_i)|\simeq |w_i|^{s_0}$.
 
Define $M_i=\max\{|z_i|^r|w_i|^s,\ (r,s)\in\mathcal{S}\}$. 
Since $|z_i|=o(|w_i|)$, the convexity
of the map $(r,s)\mapsto |z_i|^r|w_i|^s$ implies that
$z_i^rw_i^s=o(M_i)$ for $(r,s)\in\mathcal{P}\setminus\mathcal{S}$. 

We first prove that
$M_i\simeq |w_i|^{s_0}$. If $|w_i|^{s_0}=o(M_i)$, by convexity, we have $M_i\simeq 
|z_i|^{r_1}|w_i|^{s_1}$ and $z_i^rw_i^s=o(M_i)$
 for $(r,s)\in\mathcal{P}\setminus (r_1,s_1)$. 
It follows that $|f_2(p_i)|=o(M_i)=o(|f_1(p_i)|)$, which is impossible
since $f(p_i)$ is far from $X$.

Hence $M_i\simeq |w_i|^{s_0}$. Since $\Pi_2\cap\mathcal{S} =(0,s_0)$, we obtain
$|f_2(p_i)|\simeq |w_i|^{s_0}$.
\end{pf}

\begin{example} \label{ex_externalrays}
Another interesting feature, related 
to the question of laminarity, is the case where all $d_i$ 
equal 1. 
 This  holds when 
the maps have the following form
$$f(z,w)=\left(\prod_{i=1}^d(z-a_i w) + P_{d-1}(z,w), Q_{d-1}(z,w) \right),$$
where the $a_i$ are distinct and $\mathrm{deg} Q_{d-1}\leq d-1$. 
The indeterminacy set is $I= \bigcup_i
[a_i:1:0]$ and $f(L_\infty\setminus I)=[1:0:0]=X$. 
Notice this form is generic among the maps $f=(f_1,f_2)$ such that 
$\mathrm{deg}(f_1)>\mathrm{deg}(f_2)$. 
We also assume that there exist exponents $\ell_i$ such that (\ref{eq_speed}) 
holds, and that the line at infinity is attracting.

This occurs for instance when  the homogeneous
  term of highest degree of  $Q_{d-1}$ does not vanish at
  $[a_i:1]\in\pu$,  which  still 
  remains a generic property. In this case the map is regular in the sense 
  of \cite{dis2} and all $\ell_i=d-1$.

The horizontal-like maps involved 
all have degree 1, and 
for each $\alpha$ the set $\mathcal{K}_\alpha$ is a disk
 through the indeterminacy point $I_{\alpha(0)}$, transverse to $L_\infty$.
 
Our construction shows the Green current 
 is uniformly laminar  in a neighborhood of $L_\infty$
 in $\cd$, with Cantor transverse structure
 induced by the symbol space $\Sigma$. In particular the construction
  of external rays given in  \cite{bj} works in this case. 
  
 \begin{thm}\label{thm_externalrays}
 Under the preceding hypotheses, there exists a family of external rays,
 parame\-trized by $ \Sigma\times S^1$. Almost every ray lands,
 with respect to the probability
  measure $\nu\otimes d\theta$ on $ \Sigma\times S^1$,
  and if $e$ denotes the associated
 endpoint mapping, $e_*(\nu\otimes d\theta)=\mu$,
 where $\mu$ is  the maximal entropy  measure of $f$.
 \end{thm} 

\begin{pf} the partial Green function $G_\alpha$ is well
defined on $\mathcal{K}_\alpha\cap N$ and has no critical points there 
since $G_\alpha(p)= \log\norm{p}+c_\alpha+o(1)$ on the disk
$\mathcal{K}_\alpha \cap  N$ --recall $N=\bb_1\cup\cdots\cup \bb_m$--
thus the flow of the gradient vector field
$\nabla(G_\alpha\rest{\mathcal{K_\alpha}})$ determines a family of
external rays through $I_{\alpha(0)}$ indexed by $\theta\in
S^1$. 
Extending these external rays to $\mathcal{K}'$
(the part of $\mathcal{K}$ in the basin of infinity)
 involves a
cutting procedure. We only sketch the argument, see \cite[\S 6]{bj} for
details.\\

The Riemann surface $f^{-1}\mathcal{K}_{\sigma\alpha}$ has $m$
connected components near infinity, one of which is
$\mathcal{K}_\alpha$; however globally the number of components may
drop because of the branch locus. 
So we cut  $f^{-1}\mathcal{K}_{\sigma\alpha}$ into $m$ pieces by
removing the branches of external rays located behind  the critical
points of $f$. One gets $m$ simply connected Riemann surfaces
$\mathcal{K}_{\alpha,1}$ extending the $\mathcal{K}_\alpha$. 

Iterating
this procedure yields increasing sequences of simply connected
``hedgehog-like" Riemann surfaces $\mathcal{K}_{\alpha,k}$. Moreover
the partial Green function $G_\alpha$ is still well defined on  the extension
$\mathcal{K}_{\alpha,k}$ by 
$$G_\alpha= \unsur{\ell_{\alpha(0)}\cdots \ell_{\alpha(k-1)}}
G_{\sigma^k\alpha} \rond f^k.$$ One thus obtains by incresing union
a simply connected  
$\mathcal{K}_{\alpha,\infty}$ equipped with a harmonic
partial Green function $G_\alpha$. Moreover
$\mathcal{K}_{\alpha,\infty}$ is of finite area $\nu$-almost surely
because $$T\rest{\mathcal{K}'}= \int_\Sigma
[\mathcal{K}_{\alpha,\infty}] d\nu(\alpha).$$ Lemma 7.1 in \cite{bj}
asserts then that if the disk  $\mathcal{K}_{\alpha,\infty}$ has
finite area, $d\theta$-almost every external ray lands, giving then
rise to a measurable {\sl endpoint mapping} 
 $e(\alpha, \theta)$ on $\Sigma\times S^1$, and a probability measure
$e_*(\nu\otimes d\theta)$.\\

The new point here is the proof that $\mu_e=e_*(\nu\otimes d\theta)$ equals
the equilibrium measure $\mu$. Let us assume that some $\ell_i>1$, so
that $d_t>d$ and $\mu$ is well defined \cite{dis1}. 

One needs to analyze the action of $f^*$ on the space of external rays
$\Sigma\times S^1$. First, $f:\mathcal{K}_{i\alpha}\cv \mathcal{K_\alpha}$ is a
holomorphic map of degree $\ell_i$ near the point at infinity, so 
the action of $f$ on the unit tangent circle at this point is
multiplication  by $\ell_i$, and $f^*d\theta =\ell_i d\theta$. 
On  the other hand since all $d_i=1$,
$\nu$ is the balanced measure on $\Sigma$, so for the partial inverse 
of the shift $\sigma$, $\sigma_i:
i\alpha\mapsto \alpha$, one has $(\sigma_i)_* \nu=\nu$. From  this one
easily gets $f^* (\nu\otimes d\theta)  = (\sum \ell_i) (\nu\otimes
d\theta) = d_t (\nu\otimes d\theta)$.

The measure $\nu\otimes d\theta$ induces by pushing along external
rays a diffuse measure $\mu_1$
on the boundary of the first generation of disks 
$\bigcup_{\alpha\in \Sigma} \fr(\mathcal{K_\alpha}\cap N)$, and it
follows
from the foregoing discussion that $d_t^{-n}(f^n)^*\mu_1$ converges to
the endpoint measure $\mu_e$. On the other hand a $\mu_1$-generic 
point is non exceptional and $d_t^{-n}(f^n)^*\mu_1\cvf \mu$ the
equilibrium measure. This concludes the proof.
\end{pf}
  \end{example}


\section{Convergence theorems}\label{sec_conv} 

The purpose of this section is to prove Theorem \ref{thm_cv_main}
(\S \ref{subs_general})
We will also  find sufficient conditions ensuring that $\supp (T)$ is pluripolar 
(\S \ref{subs_polar}) or $T$ has continuous potential (\S \ref{subs_apriori}).
The latter estimate will allow us to derive an alternate proof of the convergence 
theorem (\S \ref{subs_largedeg}) giving explicit estimates on the rate 
of convergence, which is useful in applications.

Let us now recall the setting: 
$\set{(f_i,\bb_i,\bb_{i+1})}_{i\geq 1}$ is a sequence of meromorphic
 horizontal-like maps of respective degree $d_i$, defined in a chain of bidisks.
 As before, we always tacitly identify the abstract bidisk $\bb_i$ and $\bb$ through
 a fixed biholomorphism. 

In subsections  \ref{subs_polar} to \ref{subs_largedeg}
we will moreover assume the indeterminacy points are located in 
the distinguished
horizontal line  $(w=0)$, and the rate of attraction to this line is
controlled by exponents $\underline{\ell}_i$ and $\overline{\ell_i}$.
Of course this corresponds to the situation in section \ref{sec_speed}.

We will also suppose that the sequence of horizontal-like maps is finitely generated.
  The results may nevertheless
  be adapted to (suitable) compact families of mappings.\\


\subsection{A general convergence theorem}\label{subs_general} 

The methods in this paragraph are based on the well known fact that the slice measures
$m^w$ of a vertical positive closed current $T$ {\sl vary holomorphically} in the
following sense:
for every holomorphic function $F$ in the neighborhood of $\supp(T)$, 
the function $f$ defined by
$$f(w)=\int_{\dd\times \set{w}} F(\zeta,w)dm^w(\zeta)
= \int F(z,w) T\wedge [\dd\times \set{w}]
$$ is holomorphic. See 
e.g. Slodkowski \cite{slod} for a proof.  This can also be seen as a consequence of 
the approximation theorem for currents \cite{ds}:  $T$ can be approximated 
by a sequence of  subvarieties converging to $\supp(T)$ in the
Hausdorff topology, and for those the result is trivial. 

This provides us 
with a linear map 
$$\Phi_T: F(z,w)\longmapsto f(w)= \int_{\dd\times \set{w}} F(\zeta,w)dm^w(\zeta).$$ 
We call the {\sl moment functions} the functions $f$ obtained in this way.
We want to use holomorphic functions as test functions to prove the convergence
of sequences of currents of the form $(\el_1\cdots \el_n T_{n+1})$. Recall formula
(\ref{eq_canonical}) associates
 a canonical potential to a vertical normalized current in $\bb$.  

\begin{thm}\label{thm_cvmoments}
Let $\set{(f_i, \bb_i,\bb_{i+1})}$ be a sequence of meromorphic horizontal-like maps, 
such that $f_i^{-1}(\bb_{i+1})\subset \dd_{1-\e}\times \dd\subset \bb_i$ for  
a fixed $\e>0$. Let  $K=\bigcap\big(f_1^{-1}\cdots f_n^{-1}(\bb_{n+1})\big)$.
Assume for each $n$, $T_{n}$ is a normalized vertical current in
$\bb_{n}$. 
Then 
\begin{enumerate}[\em i.]
\item if $K$ is polynomially convex and  the currents $T_n$ have
 canonical potential uniformly  bounded with respect to $n$, 
\item or if $K$ has zero Lebesgue measure,
\end{enumerate}
then  the sequence   
$(\el_1\cdots \el_n T_{n+1})$ converges in $\bb_1$ to a current independent 
of $(T_n)$. 
\end{thm}

The hypothesis of the first item holds for instance if the maps $f_i$ are 
 holomorphic $\bb_i\cv\cd$. Moreover in that case, if for every $n$, 
 $\supp T_n \subset \bb_{n}\setminus (f^{-1}\bb_{n+1})$, the assumption on the 
potentials may be dropped (corollary \ref{cor_cvboundary}). \\

The following lemma tells us which knowledge on the moment functions is enough
to recover $T$, in case $K$ is polynomially convex or is small in some
sense.

\begin{lem}\label{lem_moments}
Let $T_1$ and $T_2$ be vertical closed positive currents in $\bb$,
supported on the vertical closed set $K$, satisfying 
$\Phi_{T_1}=\Phi_{T_2}$ --i.e. $T_1$ and $T_2$ have the same moment functions.

Assume  for $w$ but possibly in a set of zero Lebesgue measure, 
the slice $K\cap (\dd\times\set{w})$  does not separate 
$\dd\times\set{w}$, and that the slice measures are supported on $\fr(K\cap(\dd\
\times\set{w}))$. Then $T_1=T_2$.

In particular  if $K$ has  zero 
3-dimensional Hausdorff measure, and  $T_1$ and $T_2$ are supported on $K$ with
the same moment functions, then $T_1=T_2$.
\end{lem}

\begin{pf} notice first that if $T_1$ and $T_2$ are vertical closed positive 
currents such that for (Lebesgue) almost every $w$, $m^w_1=m^w_2$, then
$T_1=T_2$. Indeed the canonical 
potentials $u_1$ and $u_2$ coincide on almost every slice, hence everywhere
by subharmonicity.

Notice also that if $K$ is polynomially convex in $\bb$, so are the slices
$K\cap (\set{w}\times\dd)$, so $K$ does not
separate the slice.

The following fact is classical (see the Lebesgue-Walsh Theorem, \cite{ga} 
p.36). Let $X$ be the boundary of a polynomially convex set in $\cc$. Then 
real continuous  functions on $\fr X$ are uniform limits of harmonic polynomials.
This implies that if two real measures $m_1$ and $m_2$ are supported
on $\fr X$ and coincide on holomorphic functions, they are equal. So under the 
hypothesis of the lemma, we thus get that $m^w_1=m^w_2$ for almost every $w$. 

Assume now the Hausdorff measure $\Lambda_3(K)=0$. 
The second projection $K\cv\dd$ is surjective since
$K$ supports a vertical closed current, so 
almost every slice $K\cap (\dd\times\set{w})$ has zero length. In particular the slice
$K\cap (\dd\times\set{w})$ does not separate $\cc$ and has empty interior. 
This holds for instance if $K$ is pluripolar.  
\end{pf}

We have an analogous result in the zero measure case.

\begin{lem}\label{lem_moments2}
Let $T_1$ and $T_2$ be vertical closed positive currents in $\bb$,
with support in $K$. Assume that $K$ has zero Lebesgue measure and
that for almost every slice $\set{w}\times\dd$, and every holomorphic function
$F$ in the neighborhood of $K\cap(\set{w}\times\dd)$, one has 
$m^w_1(F)=m^w_2(F)$. Then $T_1=T_2$.
\end{lem}

\begin{pf} for Lebesgue almost every slice $\set{w}\times \dd$, 
$K\cap(\set{w}\times\dd)$ has Lebesgue measure zero. 

It is known (the Hartogs-Rosenthal Theorem, see \cite{ga})
that if  $X\subset \cc$  is closed and has zero measure, then continuous
functions on $X$ are uniform limits of rational functions with poles
off $X$. 
In particular if $m_1$ and $m_2$ are positive measures supported on
$X$ and coincide on functions holomorphic in the neighborhood of $X$,
then $m_1=m_2$. This proves the lemma.
\end{pf}

Let $\set{(f_i, \bb_i,\bb_{i+1})}$ be a (not necessarily finitely generated)
sequence of meromorphic horizon\-tal-like maps , and assume 
$f_i^{-1}(\bb_{i+1})\subset \dd_{1-\e}\times \dd\subset \bb_i$ for  
a fixed $\e>0$.
Let 
$$\mathcal{C}_n=\set{ \el_1\cdots \el_{n-1} T_{n},~T_{n}\text{ vertical normalized
in }\bb_n}$$
and for a given slice $\set{w}\times\dd$ and a holomorphic function $F$
defined in a neighborhood of $K\cap(\set{w}\times\dd)$, we
let $\mathcal{H}_n(F)$ be the range of $\Phi_T$ acting on $\mathcal{C}_n$:
$$\mathcal{H}_n(F,w)=\set{m^w(F)=\int (T\wedge[\set{w}\times \dd])F ,~T\in\mathcal{C}_n}. $$
Clearly these sequences  are nonincreasing; 
we denote by  $\mathcal{C}_\infty=\bigcap\mathcal{C}_n$ and 
$\mathcal{H}_\infty(F,w)=\bigcap\overline{\mathcal{H}_n(F,w)}$. We can now state the 
convergence result for the moment functions.

\begin{prop}\label{prop_cvmoments}
Let $\set{(f_i, \bb_i,\bb_{i+1})}$ be a sequence of meromorphic 
horizontal-like maps, 
such that $f_i^{-1}(\bb_{i+1})\subset \dd_{1-\e}\times \dd\subset \bb_i$ for  
a fixed $\e>0$.

Then with the preceding notation, for every slice $\set{w}\times\dd$
and every $F$, $\mathcal{H}_\infty(F,w)$ is a single point, that is,
all limit currents have the same moments in the slice.
\end{prop}

From lemmas \ref{lem_moments} and
 \ref{lem_moments2} one deduces the following corollaries.

\begin{cor}\label{cor_cvpolar}
If $K=\bigcap\big(f_1^{-1}\cdots f_n^{-1}(\bb_{n+1})\big)$ 
has zero Lebesgue measure (for instance if $K$ is pluripolar), 
$\mathcal{C}_\infty$ is a single current. In particular 
if for every $n$, $T_{n+1}$ is a  vertical normalized
current in $\bb_{n+1}$, the sequence $(\el_1\ldots\el_n T_{n+1})$ converges to 
the unique $\tau\in\mathcal{C}_\infty$.
\end{cor}

The pluripolarity assumption
holds in particular under the hypotheses of proposition 
\ref{prop_polar} below, i.e. if the $f_i$ have 
attracting indeterminacy points with large 
exponents $\underline{\ell}_i$.

\begin{cor}\label{cor_cvboundary}
Assume $K$ is polynomially convex and 
for every $n$, $T_n$ is a vertical normalized
current in $\bb_n$
 such that for a.e. $w$,
 $\supp(\el_1\cdots\el_nT_{n+1})\cap K\cap(\dd\times\set{w})
 =\emptyset$. Then 
the sequence $\el_1\cdots\el_nT_{n+1}$ converges.

This holds for instance if for each $n$,
$\supp T_n\subset \bb_n\setminus f_n^{-1}(\bb_{n+1}).$
\end{cor}

\noindent{\bf Proof of proposition \ref{prop_cvmoments}:}
 the proof is reminiscent of the convergence theorem in \cite{hl}.
We fix a slice $\set{w}\times\dd$ and a 
holomorphic function $F$ in a neighborhood $U$ (in $\cd$)
of $K\cap 
(\set{w}\times\dd)$. If $n_0$ is large enough, for $n\geq n_0$, 
$$\bigcap_{k\leq n}\big(f_1^{-1}\cdots f_k^{-1}(\bb_{k+1})\big) \cap
(\set{w}\times\dd) \subset U' \finc U;$$ by now we only consider
integers $n\geq n_0$. Without loss of generality, assume $\abs{F}\leq
1$ on $U'$.

For $\theta\in\dd$, let $i_\theta$ be the linear contraction 
$i_\theta(z,w)=(\theta z,w)$; $i_\theta$ acts
on vertical normalized currents by $T\mapsto (i_\theta)_*T=T_\theta$. 
If $\theta=0$, 
for any $T$, $T_0=(i_0)_*T=[z=0]$. We  let $\mathcal{H}_n=
\mathcal{H}_n(F,w)$, and 
$$\theta\cdot\mathcal{H}_n=\set{\int F (\el_1\cdots\el_n T_\theta 
\wedge[\set{w}\times \dd])
,~T \text{ vertical in }\bb_{n+1}}
\subset\mathcal{H}_n.$$
 The aim is to prove $\mathcal{H}_\infty$ is reduced to a single 
point.

 Since
each $\el_nT_{n+1}$  is a $T_{\theta_0}$ for a vertical $T$ in $\bb_n$, there exists
 a fixed $\theta_0<1$ (any $\theta_0>1-\e$ will do),
such that for
 every $n\geq 1$, 
 \begin{equation}\label{eq_theta}
 \mathcal{H}_{n+1}\subset\theta_0\cdot\mathcal{H}_n.
 \end{equation}

The crucial fact is that for fixed $n\geq n_0$ and $T$, 
$\theta\mapsto \int F (\el_1\cdots\el_n T_\theta \wedge[\set{w}\times
\dd])$ is a holomorphic function of $\theta\in \dd$, bounded by 1. 
The proof goes as follows.  The  Duval-Sibony \cite {ds} theorem allows to
approximate $T$ by vertical varieties $\unsur{m_k}[V_k]$
in $\bb_{n+1}$. Of course 
$$g_k(\zeta, \theta)=\int F \left( \unsur{m_k}
\el_1\cdots\el_n (i_\theta)_*[V_k] \wedge[\set{\zeta}\times
\dd]\right)$$ is holomorphic in the variable $(\zeta,\theta)$, for
$\zeta$ in the neighborhood of $w$ and $\theta\in \dd$, and bounded by
1 (because $F$ is bounded by 1 in $V'$). Since 
$\el_1\cdots\el_n \left(\unsur{m_k}[V_k]_\theta\right)
\cvf \el_1\cdots\el_n
T_\theta$, the convergence holds for almost every slice, and $g_k$
being holomorphic and bounded by 1, the sequence $(g_k)$ converges 
everywhere towards the desired holomorphic function.\\
 
We can thus consider $\bigcup_{\theta\in\dd}\theta\cdot\mathcal{H}_n=
\mathcal{H}_n$ as 
family of holomorphic functions in the variable $\theta$, bounded by 1, 
namely the set of 
$f_n(\theta)= \el_1\cdots\el_n T_\theta \wedge[\set{w}\times
\dd])F$ for $T$ in
$\mathbb{B}_{n+1}$ (we drop the $w$ variable). 
This family has the additional property that for fixed $n$, all 
$f_n(\theta)$ coincide for $\theta=0$, since for every vertical current
$T$, $T_0$ is the vertical line $[z=0]$.
Every $f_n\in\mathcal{H}_n$ can be embedded as
 the time 1 of a family $f_n(\theta)$.
Our first observation (\ref{eq_theta})
also yields that $f_n$ is the time $\theta_0$ of a parameter
family $f_{n-1}(\theta)\subset\mathcal{H}_{n-1}$. \\

Extracting convergent subsequences from the normal families,  relation 
(\ref{eq_theta}) persists at the limit, i.e. $\theta_0\cdot\mathcal{H}_\infty\subset 
\mathcal{H}_\infty$. More precisely for every $f\in \mathcal{H}_\infty$ there exists 
a holomorphic  $g(\theta)\subset\mathcal{H}_\infty$ with
$f=g(\theta_0)$. On the other hand every $f\in \mathcal{H}_\infty$ is the time 1 
of a holomorphic family $f(\theta)$.

For every $n$, let $f_n(\theta)\subset \mathcal{H}_n$. 
Extract a subsequence $n_j$ so that $f_{n_j}(0)$ converges. The limit $f(0)$ 
is the common value of all families $f(\theta)$ obtained as normal limits 
of $f_{n_j}(\theta)$, for $\theta=0$. Assume there exists
 a  non constant such family $f(\theta)\subset\mathcal{H}_\infty$, and let 
 $$R=\sup\set{\norm{f(\theta)-f(0)}_{L^\infty},~ f(\theta)\subset\mathcal{H}_\infty
 ,~\theta\in \dd}.$$
The Schwarz lemma asserts that for every family $f(\theta)$, 
$\norm{f(\theta)-f(0)}_{L^\infty}\leq \abs{\theta}R$. Pick  a nonconstant 
function $f$ that almost maximizes $\norm{f-f(0)}$ in
 $\mathcal{H}_\infty$, and rotate it so that $\norm{f-f(0)}=\abs{f(1)-f(0)}$
--this corresponds to replacing $T$ by some $T_\theta$ with
 $\abs{\theta}=1$.
By the Schwarz lemma, 
$f(1)$ cannot be of the form $g(\theta_0)$ and we reach a contradiction. 
We thus get that for the subsequence 
$n_j$ all normal limits of $f_{n_j}(\theta)$ are independent of
 $\theta$ (and thus of the currents $T_{n_j+1}$)
and  equal $f(0)$. Since the sequence $\mathcal{H}_n$  decreases, this remains true 
for the whole sequence $n$. \hfill $\square$\\

\noindent{\bf Proof of theorem \ref{thm_cvmoments}:} under the
 hypothesis that $K$ is polynomially convex, 
we split the proof into two cases.

Assume first the sequence $(d_1\cdots d_n)_{n\geq 1}$ is bounded. This means 
that for large enough $n$, $d_n=1$, so 
we can assume all  $d_n=1$. In this case the theorem is a version of the 
 Stable  Manifold Theorem. The proof goes as follows. Starting with a vertical 
holomorphic
graph $\Gamma_{n+1}$ in $\bb_{n+1}$, $\el_1\cdots\el_n[\Gamma_{n+1}]$ 
is also a vertical 
graph, as well as the cluster values of this sequence. Since a 
vertical graph is determined 
by its moment functions, from proposition \ref{prop_cvmoments} we get 
$\el_1\cdots\el_n[\Gamma_{n+1}]\cv [\Gamma_\infty]$. Now it is clear that 
 the limit set $K=\bigcap  f^{-1}_1\cdots 
f^{-1}_n(\bb_{n+1})\subset \Gamma_\infty $ so the convergence holds
for every 
starting
 current $T_{n+1}$, and $K=\Gamma_\infty$.  \\
 
 Suppose now $d_1\cdots d_n\cv\infty$. We will prove that for any cluster point 
 $\tau$ of the sequence 
 $\el_1\cdots \el_{n} T_{n+1}$, its slice measures on $\dd\times \set{w}$
 are supported on $\fr(K\cap (\dd\times\set{w}))$,
for a.e. $w$. The polynomial convexity of $K\cap
 (\dd\times\set{w})$ is clear. 
 By lemma \ref{lem_moments} this 
 implies $\tau$ is uniquely defined by its moment functions,  
thus by proposition \ref{prop_cvmoments}, $(\el_1\cdots \el_n 
 T_{n+1})$ converges.\\
 
If $\el_1\cdots \el_{n_j} T_{n_j+1}\cvf\tau$,
by Slicing theory, 
for almost every $w$, the sequence of slice measures   
$\el_1\cdots \el_n T_{n_j+1}\wedge [\dd\times\set{w}]$ converges
to $\tau\wedge  [\dd\times\set{w}]$. We pick up such a $w$.
Assume  the $T_n$ have uniformly bounded potentials, and let $V$ be an open 
subset in the interior of $K\cap (\dd\times\set{w})$.
 For every $n_j$, $f_{n_j}\rond \cdots\rond f_1$ is  a well defined bounded family 
of holomorphic mappings $V\cv \bb_{n_j+1}$ --use extension across indeterminacy 
points if any. Let $u_{n_j+1}$ be the canonical potential of $T_{n_j+1}$, 
then 
$$\unsur{d_1\cdots d_{n_j}}  u_{n_j+1}\rond f_{n_j}\rond \cdots\rond f_1 $$ 
is a potential of the slice  of $\el_1\cdots \el_{n_j}
 T_{n_j+1}\wedge [\dd\times\set{w}]$ 
in $V$. Since $u_{n_j+1}$ is bounded this sequence uniformly converges to zero 
 and the slice of $\tau$ gives zero mass to $V$.\\
 

The second item of the theorem 
is corollary \ref{cor_cvpolar}. 
   \hfill $\square$\\

\begin{rmk}\normalfont \label{rmk_Pesin} 
When all $d_i=1$, there is no need to use proposition \ref{prop_cvmoments}.
Applying the Schwarz lemma in the unit ball
of $H^\infty(\dd)$ viewed as the space of vertical graphs 
is enough. This also works in higher dimension.
This  provides  a particularly simple proof
of the  Stable Manifold Theorem for 
hyperbolic periodic  points of holomorphic mappings; similarly,  allowing
the source and target polydisks to differ at each step gives  the Pesin 
Stable Manifold Theorem. 
\end{rmk}


\subsection{Pluripolarity of $K$ in the case of large exponents}\label{subs_polar}

In this paragraph we prove
that if the maps  have indeterminacy points  and 
the Lojaciewicz exponents $\underline{\ell}_i$ dominate the
$d_i$, then
the  set $\bigcap  f^{-1}_1\cdots 
f^{-1}_n(\bb_{n+1})$ must be pluripolar. This forces the limiting potentials to
be singular.  The hypotheses in this paragraph are designed to fit 
with the maps considered in section \ref{sec_speed}.

We assume the maps $f_i$ have exactly one
indeterminacy point at the origin, and 
moreover satisfy
\begin{equation}\label{eq_indet}
f_i(0)\cap\bb_{i+1}\neq\emptyset\text{ and }((w=0)\backslash\set{0})\cap f_i^{-1}
\bb_{i+1}=\emptyset
\end{equation}
As a consequence $f_i(0)\cap\bb_{i+1}$ is a 
non trivial one dimensional
subvariety, which is necessarily horizontal due to item ii. of the
definition of horizontal-like maps.   

We also make for simplicity the assumption that the sequence $(f_i)$ 
has  finitely many generators.
This implies uniformity with respect to $i$ of the constants
appearing in the estimates.

We use the following notations:
$K_n=  f^{-1}_1\cdots 
 f^{-1}_{n}(\bb_{n+1})\cup\set{0}$ and  $K=\bigcap_n K_n$; both are vertical subsets. 
One easily proves that under the preceding hypotheses
$\overline{K}_{n+1}\subset K_n$, so
$K$ is closed. 
 
\begin{prop}\label{prop_polar}
Let  $\set{(f_i,\bb_i,\bb_{i+1})}_{i\geq 1}$ 
 be a finitely generated
  sequence of meromorphic horizontal-like maps with
  exactly one indeterminacy point at 0, satisfying  condition (\ref{eq_indet}). 
  We assume the line $(w=0)$ is attracting in the following sense 
  \begin{equation}\label{eq_attract2}
  \forall i,~ \forall p=(z,w)\in f_i^{-1}(\bb_{i+1}),~   
\dist(f_i(p), (w=0))  \leq \abs{w}^{
  \underline{\ell}_i}.
  \end{equation}
Then if $$\limsup_{n\cv\infty} 
\frac{\underline{\ell}_1\cdots \underline{\ell}_n}{d_1\cdots
 d_n}=\infty,$$ the limit set $K$ is pluripolar.
\end{prop}

The $\abs{w}^{\underline{\ell}_i}$ term on the right hand side of 
condition (\ref{eq_attract2}), may be replaced by $C_i \abs{w}^{ \underline{\ell}_i}$
provided all $ \underline{\ell}_i\geq 2$. In the general case  where some 
$\underline{\ell}_i=1$ we need to assume  $C_i\leq 1$.

We begin with a couple of lemmas.

\begin{lem}\label{lem_cone}
Let $(f,\bb_1, \bb_2)$ be a meromorphic horizontal-like map with
  exactly one indeterminacy point at 0, and satisfying
  condition (\ref{eq_indet}).
  Then there exists a
constant $0<\tau\leq 1$ such that for every $p=(z,w)\in f^{-1}(\bb_2)$
  close enough to 0,
  $\norm{p}\leq \abs{w}^\tau$.
  \end{lem}

\begin{pf} we prove there exists an open analytic cone $\mathcal{C}$
 of equation $\abs{w}<\abs{z}^n$ escaping under $f$, that is 
 $f(\mathcal{C})\cap\bb_2=\emptyset$. This means the non escaping points 
 satisfy the inequality $\abs{w}\geq \abs{z}^n$ which implies the lemma.
 
 Recall the meromorphic map $f$ writes as $\pi_2\rond\pi_1^{-1}$ where $\pi_1$
 is a composition of point blow-ups. Let $E$ be the proper transform of the 
 line $(w=0)$ by $\pi_1$. The second item of
condition (\ref{eq_indet}) means $\pi_2$ sends $E$ outside
 $\bb_2$. Since $\pi_2$ is holomorphic this is true for some neighborhood 
 $N(E)$ of  $E$ in $\pi_1^{-1}(\bb_1)$. We need to
 prove the existence of  $N(E)$ such that
 $\pi_1(N(E))$ has the desired equation in the neighborhood of $0\in \bb_1$. 
 
 The exceptional divisor $\pi_1^{-1}(0)=E_1$ is a chain of rational curves intersecting
 $E$ at some point $q_0$. In the neighborhood of $q_0$, $\pi_1$ is a composition
 of  blow ups , with exceptional divisor $E_1$ intersecting $E$
 transversely.  So there 
 exists  a system of coordinates $(u,v)$ in the neighborhood of $q_0$, such that 
 $E_1=(u=0)$, $E=(v=0)$ and
 $\pi_1(u,v)=(u^a(1+o(1)),u^bv^c(1+o(1)))$. 
 In this case the image of a small 
 bidisk around $q_0$ has the desired form in $\bb_1$.
 \end{pf}
 
\begin{lem}\label{lem_pushpull}
 Let $(f,\bb_1,\bb_2)$ 
 be a meromorphic horizontal-like map with
  exactly one indeterminacy point at 0, and satisfying
  condition (\ref{eq_indet}). Then for every
  normalized vertical current $T$ in $\bb_{2}$, the current $\el T$ has Lelong number
  $\geq c(f)>0$ at the origin.
\end{lem}

\begin{pf} actually we show that  any normalized vertical closed
  positive current in $f^{-1}(\bb_2)$ has positive Lelong number at 0.
  Indeed assume $T$ is
  such a current, and let $u$ be its logarithmic potential. Then
$$u(z,w)=  \int_{\dd\times\set{w}} \log\abs{z-\zeta} dm^w(\zeta) \leq
 \log \max\set{\abs{z-\zeta},~ \zeta \in \supp(m^w)} \leq
\log\left( \abs{z}+\abs{w}^\tau\right)$$ where the final estimate follows from the
 preceding lemma and the triangle inequality. 

In particular one has $u(p)\leq \tau\log\norm{p}+O(1)$ for
$p=(z,w)$ near 0, which implies the lemma.
\end{pf}

\noindent{\bf Proof of proposition \ref{prop_polar}:} 
 consider a vertical
  normalized current $T_{n+2}$ in $\bb_{n+2}$, then by the preceding lemma
  $\el_{n+1}T_{n+2}=T_{n+1}$ is a normalized vertical
  current in $\bb_{n+1}$ with positive Lelong number at the origin. The current
  $\el_1\cdots\el_n T_{n+1}$ has support in $K_{n+1}\subset K_n\subset \bb_1$.
  
  There are two natural potentials for this current: first the
  logarithmic potential $u_{1,n}$
  defined as in (\ref{eq_canonical}), next the potential, 
which is only defined in $K_n$
  $$u_{2,n}=\frac{u({T_{n+1}})\rond f_n\rond\cdots\rond f_1}{d_1\cdots d_n},$$ where
$u({T_{n+1}})$ is the logarithmic potential  of $T_{n+1}$. 
If $p\in  \fr K_n$, then $f_n\rond\cdots\rond f_1(p)\in\fr_v
\bb_{n+1}$, so  $u_{2,n}(p)\geq (\log\e)/(d_1\cdots d_n)$. 
On the other hand $u_{1,n}\leq \log 2$
everywhere 
in $\bb_1$. This implies that $u_{1,n}\leq u_{2,n}
+O(1)$ on $ \fr K_n$. Since $u_{1,n}-u_{2,n}$ is
pluriharmonic, this inequality extends to $K_n$.

Now if $p=(z,w)\in K_{n+1}$, the second coordinate of 
$f_n\rond\cdots\rond f_1(p)$ is smaller than 
$\abs{w}^{\underline{\ell}_1\cdots \underline{\ell}_n}$. By the preceding lemma,
there exists a constant $\tau$ such that for every 
$q\in f^{-1}_{n+1}(\bb_{n+2})$, $\norm{q}\leq \dist(q,(w=0))^\tau$, 
 so we get 
$\norm{f_n\rond\cdots\rond f_1(p)}\leq 
\abs{w}^{\tau\underline{\ell}_1\cdots \underline{\ell}_n}$, and since $u({T_{n+1}})$ 
has a  positive Lelong number at 0, for $q$ near 0 one has $u({T_{n+1}})
(q)\leq c\log\norm{q} + O(1)$. Hence
\begin{equation}\label{eq_polar}
u_{1,n}(z,w) \leq u_{2,n}(z,w)+O(1)
\leq C\frac{\underline{\ell}_1\cdots \underline{\ell}_n}{d_1\cdots
 d_n} \log{\abs{w}} + O(1)
\end{equation}
 on $K_{n+1}$.
 
Extract first a subsequence $n_i$
 so that $u_{1,n_i}\leq -i^3$ on $K_{n_i+1}$, and let 
 $$u=\sum_{i=1}^\infty \unsur{i^2} u_{1,n_i}.$$
 Then $u$ is p.s.h. and $u=-\infty$ on $K$. On the other hand, since for every 
 $i$, $u_{1,n_i}\geq (\log\e)/(d_1\cdots d_n)$ 
near $\fr_v\bb_1$, $u$ is uniformly bounded there. Hence
 $K$ is pluripolar.  
 \hfill $\square$

\subsection{A  priori bounds for potentials}\label{subs_apriori}

Here we give  bounds on the potentials of the  currents 
$\el_1\cdots \el_n T$,
when the rates of attraction $\ell_i$ of indeterminacy points are controlled
by the $d_i$. This also covers the case where 
there are no indeterminacy points.
This will allow us in the next paragraph
 to give another approach to the convergence result.

The results can easily be
transposed in the setting of sequences  of polynomial-like maps in $\cc$.
As far as we know, the result is new even in that case. It gives a
lower estimate on the size of the ``filled Julia set" associated to such a
sequence, in a potential theoretic sense. 

If $T$ is a vertical (closed positive)
current in $\bb$, we let $u(T)$ be its canonical
potential as given by equation (\ref{eq_canonical}) (coordinates are 
swapped for vertical currents). 

\begin{lem}\label{lem_apriori1}
Suppose $T$ is a normalized
 closed positive current with support in
$\dd_{1-\e}\times\dd\subset \bb_{i+1}$, and $f_i$ is a meromorphic
horizontal-like map with  possible indeterminacy points in the line $(w=0)$.
 Then for every horizontal slice $\dd\times
\set{w}$ in $\bb_i$, $w\neq 0$,  one has
$$\abs{ u(\el_i T) -\unsur{d_i} u(T)\rond f_i}_{(\dd\times
\set{w})\cap f_i^{-1}(\bb_{i+1})} \leq C_i(1+e_i\abs{ \log\abs{w}}),$$ 
where $C_i$ and $e_i$ depend only on $f_i$ and $\e$ ($e_i=0$ if $f_i$ is 
holomorphic).
\end{lem}

\begin{pf} note first that the function 
$ u(\el_i T) -\unsur{d_i} u(T)\rond f_i$ is pluriharmonic where it is
defined, i.e. in $f_i^{-1}(\bb_{i+1})$, so to estimate it on a given
horizontal slice, it suffices to  restrict to the boundary
$(\dd\times\set{w})\cap\fr \left(f_i^{-1}(\bb_{i+1})\right)$.

Now if $p\in \fr \left(f_i^{-1}(\bb_{i+1})\right)$, $f_i(p)\in
\fr_v\bb_{i+1}$, so $\dist(f_i(p),\supp (T))\geq \e$, and $u(T)(f_i(p))\geq
\log \e$, which implies 
\begin{equation}\label{eq_lem1}
\abs{\unsur{d_i} u(T)\rond f_i}_{(\dd\times
\set{w})\cap \fr f_i^{-1}(\bb_{i+1})}\leq \unsur{d_i}\abs{\log \e}
\end{equation}
On the other hand, 
\begin{equation}\label{eq_lem2}
\abs{u(\el_i T)}_{ (\dd\times\set{w})\cap\fr(f_i^{-1}(\bb_{i+1})) }
\leq \abs{ \log\left( 
{\rm horizontal ~dist.}\left(\supp \el_i T, \fr f_i^{-1}(\bb_{i+1})~\right)
\right)},
\end{equation}
while the horizontal distance is evaluated using the mean value
inequality. Indeed, if $p_1\in\supp \el_i T$ and $p_2\in  \fr
f_i^{-1}(\bb_{i+1})$ are on the same horizontal line 
$$\e\leq \abs{f_i(p_1)- f_i(p_2)} \leq \frac{C_1}{\abs{w}^{e_i}}
  \abs{p_1-p_2},$$ for some constants $C_1$ and $e_i$.
Hence $\abs{p_1-p_2}\geq C_2\abs{w}^{e_i}$, which implies  
$$\abs{u(\el_i T)}_{ (\dd\times\set{w})\cap\fr(f_i^{-1}(\bb_{i+1}))}\leq C(1
+e_i\abs{ \log\abs{w}})$$ by (\ref{eq_lem2}).
This, together with (\ref{eq_lem1}) yields the desired estimate.
\end{pf}

We now iterate the preceding lemma for a sequence of mappings
 $\set{(f_i,\bb_i,\bb_{i+1})}_{i\geq 1}$. Remember  the constants $C_i$
 and $e_i$ are uniformly bounded, for we consider a finitely generated 
  family $(f_i)$. 

\begin{lem}\label{lem_apriori2}
If $T$ is a normalized
vertical current in $\dd_{1-\e}\times\dd\subset \bb_{n+1}$, 
 one has for every $p\in f^{-1}_1\rond\cdots 
\rond f^{-1}_n(\bb_{n+1})$,
$$\abs{u(\el_1\cdots\el_n T)(p)} \leq C
 \sum_{i=0}^{n-1} \frac{1+ e 
  \big|\log{\abs{(f_i\rond \cdots \rond f_1(p))_2}} \big| }{d_1\cdots
  d_i} + \frac{u(T)(f_n\rond \cdots \rond f_1(p))}{d_1\cdots d_n}
,$$ where the subscript $2$ means only the log of the second coordinate is
involved, and $e=0$ if the maps are holomorphic.
\end{lem}

\begin{pf} the proof is just applying inductively the preceding lemma, 
by using the following trick: in $f^{-1}_1\cdots 
f^{-1}_n(\bb_{n+1})$, one has 
$$u(\el_1\cdots\el_n T) = u(\el_1\cdots\el_n T)- 
\frac{u(T)\rond f_n\rond \cdots \rond f_1}{d_1\cdots d_n} + 
 \frac{u(T)\rond f_n\rond \cdots \rond f_1}{d_1\cdots d_n}$$ which
 splits as 
$$ \sum_{i=1}^{n}\left[
\frac{u(\el_{i}\cdots\el_n T)\rond f_{i-1}\rond \cdots \rond f_1}{d_1\cdots d_{i-1}}- 
\frac{u(\el_{i+1}\cdots\el_n T)\rond f_i\rond \cdots \rond
  f_{1}}{d_1\cdots d_{i}}\right]
   + \frac{u(T)\rond f_n\rond \cdots \rond f_1}{d_1\cdots d_n}$$
with an obvious meaning when $i=1$ or $i=n$, i.e. $f_0=id$, $d_0=1$,
etc.  This rewrites as
\begin{equation}\label{eq_atrroce}
\sum_{i=1}^{n} \unsur{d_1\cdots
  d_{i-1}}\left[u(\el_{i}\cdots\el_n T)-\unsur{d_i}
u(\el_{i+1}\cdots\el_n T)\rond f_i \right] \rond f_{i-1}\rond
\cdots \rond f_1
 + \frac{u(T)\rond f_n\rond \cdots \rond f_1}{d_1\cdots d_n}. 
\end{equation}
Applying  lemma \ref{lem_apriori1} to the term between brackets gives the
desired bound. 
\end{pf}

From this lemma one easily deduces a priori estimates for the
sequences of potentials.
 In presence of indeterminacy points, one cannot hope for
uniform estimates on the potentials,  nevertheless we obtain
bounds on horizontal slices,
in case  the iterates do not approach the indeterminacy locus too fast.
Our approach actually also yields continuity of the limiting potentials.

\begin{prop}\label{prop_apriori}
Assume $\set{(f_i,\bb_i,\bb_{i+1})}_{i\geq 1}$ is a finitely generated sequence
of meromorphic horizontal-like maps. We moreover assume
\begin{itemize}
\item[-] either the maps have no indeterminacy points
\item[-] or  the possible indeterminacy locus is contained in 
$(w=0)$, and  the following is true for every $i$:
\begin{equation}\label{eq_attract}
\exists C_i>0,~\exists \overline{\ell}_i\geq1,~ 
\forall p=(z,w)\in f_i^{-1}(\bb_{i+1}),~ C_i^{-1}\abs{w}^{
  \overline{\ell}_i}\leq 
  \mathrm{dist}(f_i(p), (w=0)).
 \end{equation}
 \end{itemize}
Then if
$T_{n+1}$ is a vertical normalized current in $\bb_{n+1}$ with uniformly
 bounded potential,
there exist constants $C$ and $e$ (with $e=0$ when the maps are
holomorphic) such that on every horizontal slice $\dd\times\set{w}$,
$w\neq 0$, one has 
\begin{equation}\label{eq_estim}
\big| u(\el_1\cdots\el_n T_{n+1}) \big|_{\dd\times\set{w}} \leq
C\sum_{i=0}^n \frac{1+ e
\overline{\ell}_1\overline{\ell}_2\cdots \overline{\ell}_i \left(i+
\abs{\log\abs{w}} \right)}{d_1\cdots d_i}.
\end{equation}
In particular if the series on the right hand side converges, 
any cluster value $\tau$ of  $(\el_1\cdots\el_n T_{n+1})$ has continuous
 potential outside $(w=0)$.
\end{prop}

We will prove in theorem \ref{thm_cvlarge} below that under the
hypotheses of the proposition, the sequence
 $(\el_1\cdots\el_n T_{n+1})$ converges. 
The $(i+ \abs{\log\abs{w}})$  factor  may be replaced by $1+\abs{\log\abs{w}}$ 
in many cases, see  the
proof below. We stress condition (\ref{eq_attract}) is of no importance 
when there are no 
indeterminacy points.\\

\begin{pf} we prove (\ref{eq_estim}) first.
The potential of the normalized current $u(\el_1\cdots\el_n T_{n+1})$ 
is subharmonic on every horizontal slice $\dd\times\set{w}$, and its
laplacian has support in  $f^{-1}_1\cdots 
 f^{-1}_n(\bb_{n+1})$, so it suffices to control the $L^\infty$ norm
there.

We  apply lemma \ref{lem_apriori2}.
 If $p=(z,w)\in  f^{-1}_1\cdots 
 f^{-1}_n(\bb_{n+1})$ we have to estimate the second coordinate of
 $ \log{\abs{(f_i\rond
      \cdots \rond f_1(p))}}$ using the left hand side of 
equation (\ref{eq_attract}), where the constant $C_i$ is replaced by
a uniform constant. An easy induction shows that 
$$  \big|\log{\abs{(f_i\rond
      \cdots \rond f_1(p))_2}}~ \big| \leq  C \left(\sum_{j=1}^{i+1}
      \overline{\ell}_j\cdots \overline{\ell}_i \right) +   
 \overline{\ell}_1\cdots \overline{\ell}_i~ \big|\log\abs{w}\big|.$$

In case  $\overline{\ell}_i \geq 2$ the constant term can not become
overwhelming because of the simple inequality 
$$1+ \overline{\ell}_i+
\overline{\ell}_{i-1}\overline{\ell}_i
+\cdots  +\overline{\ell}_1\overline{\ell}_2\cdots \overline{\ell}_i 
\leq 2\overline{\ell}_1\overline{\ell}_2\cdots \overline{\ell}_i.$$ In
the general case
$\ell_i\geq 1$ 
one  uses the following
$$1+ \overline{\ell}_i+\cdots  +\overline{\ell}_1\overline{\ell}_2\cdots \overline{\ell}_i 
\leq  (i+1) \overline{\ell}_1\overline{\ell}_2\cdots \overline{\ell}_i.$$
Notice that in general the growth of the products
$\overline{\ell}_1\overline{\ell}_2\cdots \overline{\ell}_i$  is exponential
so  the linear term  $(i+1)$ does not affect it too much.\\

Let $\tau$ be a cluster value of the sequence  
$(\el_1\cdots\el_n T_{n+1})$ in $\bb_1$. If the 
series in (\ref{eq_estim}) converges,  $\tau$ has locally bounded potential 
$u(\tau)$ outside $(w=0)$. Starting with e.g. smooth $T_n$,
we prove $u(\tau)$ is continuous. 
We need only prove the continuity of $u(\tau)$ at $p\in K$. Let $\e>0$, and 
fix $m$ such that the remainder of the series in (\ref{eq_estim}) is smaller than $\e$. 
If $q$ is close enough to $p$, then $q\in f_1^{-1}\cdots f_n^{-1}(\bb_{n+1})$ for 
some $n\geq m$. Rewrite (\ref{eq_atrroce}) as 
$u(\el_1\cdots\el_n T)=S_m+R_m$ where $S_m$ consists 
of the $m$ first terms of the sum, so that
$\abs{R_m(q)}< \e$ and $\abs{R_m(p)}<\e$.
To get equicontinuity of the sequence $u(\el_1\cdots\el_n T)$ near $p$ 
it remains to bound $\abs{S_m(p)-S_m(q)}$. 

Only the term between brackets 
in $S_m$ depends on $n$. Let
$\tilde{T}=\el_{i+1}\cdots\el_nT_{n+1}$, the function 
$u(\el_i\tilde{T})-\unsur{d_i} u(\tilde{T})\rond f_i$ 
is pluriharmonic in $f_i^{-1}(\bb_{i+1})$ and locally  bounded,
independently on $n$ (see lemma \ref{lem_apriori1}). So the modulus of 
continuity of this function in
  $f_i^{-1}\cdots f_n^{-1}(\bb_{n+1})$ does not depend on $n$.
This concludes the proof. 
\end{pf}

\subsection{An alternate approach to the convergence theorem} \label{subs_largedeg}

We keep in this paragraph the hypotheses of Proposition
\ref{prop_apriori}. We prove a convergence theorem in the spirit of 
theorem \ref{thm_cvmoments}, with however hypotheses of different nature. 

\begin{thm}\label{thm_cvlarge}
Under the  hypotheses of Proposition
\ref{prop_apriori}, assume for every $n$, $T_{n+1}$ 
is a  normalized vertical positive closed current in $\bb_{n+1}$ 
with uniformly bounded canonical 
potential. Then if the series 
\begin{equation}\label{eq_series}
\sum_{n\geq 1}\frac{1+en\overline{\ell}_1\cdots
  \overline{\ell}_n}{d_1\cdots d_n}
\text{ (with } e=0 \text{ if the maps are holomorphic)}
\end{equation}
 converges,
  so does the
  sequence of currents $(\el_1\cdots\el_n T_{n+1})$. 
 The limit
  is independent of the currents $T_n$, and has continuous potential
 outside $(w=0)$.
\end{thm}   

As opposite to theorem \ref{thm_cvmoments}, this approach
provides explicit bounds on the rate of convergence, which is important 
for applications.  This is new even in the case of one single iterated 
H{\'e}non-like map, and   solves the question raised
 in \cite[Remark 4.4]{hl}.

 The proof yields 
  an estimate on the rate of convergence in terms of the expression
\begin{equation}\label{eq_atroce}
\unsur{d_1\cdots d_{m}}\sum_{i=1}^{\infty}\frac{1+e\overline{\ell}_{m+1}\cdots
  \overline{\ell}_{m+i}(i + m\overline{\ell}_1\cdots
  \overline{\ell}_{m})}{d_{m+1}\cdots d_{m+i}}.
\end{equation}
  Actually, the convergence statement follows from the
 assumption that (\ref{eq_atroce}) tends to zero when
  $m\cv\infty$. As before (see the proof of proposition \ref{prop_apriori}),
  the assumptions can be slightly weakened when all 
  $\overline{\ell}_i\geq 2$.\\

\begin{pf} the
letter $C$ denotes a ``constant" that may change from line to line,
independently of $m$ and $n$.
Let $n>m$
and $T_1$ and $T'_1$ be respectively vertical normalized currents in
$\bb_{m+1}$ and $\bb_{n+1}$ with uniformly 
bounded potentials. Let $S$ be a
smooth positive closed current, with $\supp S\cap (w=0)=\emptyset$, 
and $\varphi$ a real valued test function
in $\bb_1$. Recall a vertical positive closed $T$ is uniquely  determined by
the values of $\langle T, \varphi S\rangle$, for such $S$ and
$\varphi$, since these determine the slice measures outside $(w=0)$.

We need to prove $$\langle \el_1\cdots \el_m T_1 - \el_1\cdots
\el_n T'_1, \varphi S\rangle = \langle \el_1\cdots \el_m (T_1-T_2), \varphi S\rangle$$
 is small when $m$ is large (Cauchy
criterion), with $T_2= \el_{m+1}\cdots \el_n T'_1$.
By proposition \ref{prop_horiz}, $T_1-T_2= d\alpha$ for a vertically
supported $\alpha$, and we can write  
$$\left\langle \el_1\cdots \el_m (T_1-T_2), 
\varphi S\right\rangle =\left\langle\el_1\cdots \el_m d\alpha, 
\varphi S\right\rangle = \left\langle \el_1\cdots \el_m
 \alpha , d\varphi\wedge S\right\rangle .$$
Next, we split $\alpha$ as $ \alpha^{1,0}+\alpha^{0,1}$ and 
 use the Schwarz inequality
\begin{align}\label{eq_align}
\abs{ \Bigl\langle  \el_1\cdots \el_m \alpha^{1,0} , d\varphi\wedge S\Bigl\rangle} &\leq 
\Bigl\langle S, i\fr \varphi\wedge \overline{\fr} \varphi \Bigl\rangle
^{\unsur{2}}
\left\langle S, i \el_1\cdots \el_m (\alpha^{1,0}) \wedge  \el_1\cdots \el_m 
(\overline{\alpha^{1,0}})\right\rangle^{\unsur{2}}\\ \notag
&= \unsur{(d_1\cdots d_m)^{\frac{1}{2}} }
\Bigl\langle S, i\fr \varphi\wedge \overline{\fr} \varphi\Bigl\rangle
^{\unsur{2}}   \left\langle  \tel_m\cdots \tel_1 S, i
 \alpha^{1,0} \wedge\overline{\alpha^{1,0}}\right\rangle^{\unsur{2}},
\end{align} 
 where $\tel_m\cdots \tel_1 S$ is a horizontal current in $\bb_{m+1}$.
It remains to control $\langle  \tel_m\cdots \tel_1 S, i
 \alpha^{1,0} \wedge\overline{\alpha^{1,0}}\rangle$: this is the role
 of  proposition \ref{prop_apriori}, together with the following
  lemma. 

\begin{lem}[Polarization inequalities]\label{lem_polarization}
If $T_i$, $i=1,2$, (resp. $S'$) 
are vertical (resp. horizontal)
normalized positive closed currents, with 
$\supp(T_i)\subset \dd_{1-\e}\times \dd$. 
Write $T_i= dd^cu_i$, $i=1,2$ where $u_i$ is the canonical potential, and 
assume the $u_i$ are bounded on $\supp S'$.
If $\alpha=  \chi (d^cu_1-d^cu_2) + \theta d\chi$ is
as in proposition \ref{prop_horiz}, then there exists a constant $C_1$
depending only on $\e$ such that 
$$\left\langle S', i\alpha^{1,0}
    \wedge\overline{\alpha^{1,0}}\right\rangle
\leq C_1 \max\left( \norm{u_1}_{L^\infty(\supp S')},
  \norm{u_2}_{L^\infty(\supp  S')},1\right)  .$$
\end{lem}

We use the lemma to finish  the proof of the theorem. 
$T_1$ has bounded potential and
proposition \ref{prop_apriori} asserts that on the horizontal slice
$\dd\times\set{w}$
 in $\bb_{m+1}$,  
the potential of $T_2$ is controlled by 
$$C\sum_{i=0}^{n-m} \frac{1+ e
\overline{\ell}_{m+1}\cdots \overline{\ell}_{m+i} \left(i+
\abs{\log\abs{w}} \right)}{d_{m+1}\cdots d_{m+i}}.$$ 
 Lemma \ref{lem_polarization} requires a bound for the potentials
 on $\supp(\tel_m\cdots\tel_1 S)$, so here again we use condition
 (\ref{eq_attract}). 
More precisely if $\dist (\supp S, (w=0))=r$, 
one gets as in the proof of proposition \ref{prop_apriori} 
\begin{align*}
\log \dist( \supp(\tel_m\cdots\tel_1 S),  (w=0)) &\geq
-C(\overline{\ell}_m+\cdots +\overline{\ell}_m\cdots \overline{\ell}_2
 ) +\overline{\ell}_m\cdots 
  \overline{\ell}_1 \log r\\ 
 &\geq - C \overline{\ell}_m\cdots \overline{\ell}_2
  \overline{\ell}_1 ( m + \abs{\log r}). 
\end{align*}
 The current $S$ being
  fixed, this term is $O(  m\overline{\ell}_1
  \cdots\overline{\ell}_m)$. 
By (\ref{eq_align}), this gives a bound on 
$\langle  \el_1\cdots \el_m \alpha^{1,0} , d\varphi\wedge S\rangle$ in
  terms of  the 
square root of the term in (\ref{eq_atroce}). Thus if the series in 
(\ref{eq_series}) converges, Cauchy's criterion is satisfied, which
  concludes the proof.
\end{pf}

\noindent{\bf Proof of lemma \ref{lem_polarization}:} this lemma is a
version of the classical Chern-Levine-Nirenberg inequalities. Notice
however that the inequality is sharper than expected --namely the left
hand side is usually squared--
 because of the specific geometric situation. 
Assume 
$u_1$ and $u_2$ are bounded on $\supp S'$.
Since 
$\alpha=  \chi (d^cu_1-d^cu_2) + \theta d\chi$, the (1,0) part of
$\alpha$ is given by $ \chi (\fr u_1-\fr u_2) +
\theta \fr\chi$. 
So $\langle S', i\alpha^{1,0}
    \wedge\overline{\alpha^{1,0}}\rangle$ contains terms of the form 
$$\langle \chi^2 S', i\fr u_i\wedge \overline\fr u_j\rangle,  ~
\langle \chi S', i\fr u_i\wedge \theta\overline\fr\chi \rangle\text{
  and }\langle S', \theta^2\fr \chi \wedge\overline\fr\chi\rangle.$$
The form $\fr\chi$ has support in  $(\dd_{1-\e/2}\backslash
\dd_{1-\e})\times\dd$, and $\theta$ and $\fr u_i$ are controlled in
$L^\infty$ norm there, so, $S'$ being normalized,
 the terms of the second and third type are
bounded independently of $T_1$ and $T_2$. \\

The remaining issue is to estimate terms of the form 
$\langle \chi^2 S', i\fr u_i\wedge \overline\fr u_j\rangle$,
with $i,j=1,2$. Note that
using the Schwarz inequality with respect to the positive current
$\chi^2S'$, reduces to expressions like 
$\langle \chi^2 S', i\fr u_i\wedge \overline\fr u_i\rangle$. Write $u=u_i$. 
We use the classical polarization identity 
$$2i\fr u\wedge \overline\fr u = i\fr\overline\fr (u^2) - 2 u i\fr\overline\fr u$$
(notice that $u$ is not globally bounded, so we cannot assume $u^2$ is
psh by replacing $u$ by $u+C$). 

We infer (using $0\leq \chi\leq1$)
$$\abs{\langle \chi^2 S', u i\fr\overline\fr u\rangle}
\leq  \norm{u}_{L^\infty(\supp S')}\langle S',  i\fr\overline\fr
u\rangle
=  \norm{u}_{L^\infty(\supp S')},$$
since $S$ and $i\fr\overline\fr u$ are normalized, so $\int S'\wedge
i\fr\overline\fr u=1$. The other term is
estimated as follows 
$$\abs{\langle \chi^2 S', i\fr\overline\fr (u^2)\rangle}\leq  
\norm{u}_{L^\infty(\supp S')\cap\supp(\fr\overline\fr\chi)}^2 
\abs{\langle  S', i\fr\overline\fr
  (\chi^2)\rangle}\leq C.$$
This concludes the proof.\hfill $\square$

\section{Miscellaneous applications}\label{sec_appl}

We describe a couple of situations involving sequences of horizontal-like
maps, where the results of the preceding section may prove useful.
  
\subsection{Newton's method.}
Here we discuss a series of results by Y. Yamagishi \cite{ya}. 
Applying Newton's 
method to find the common zeros of a pair of polynomials $(P_1,P_2)$ yields 
a superattractive cycle when the common root is simple, and an indeterminacy
point when multiple. The dynamics near such a point looks very much like the situation
encountered in section \ref{sec_speed}.

Yamagishi \cite{ya} treats the case of a multiple root at 0 of the form
$$f= (P_1,P_2)= (z+O(\norm{p}^2), w^2-z^2+ O(\norm{p}^3))\text{, where } p=(z,w);$$
he proves that  the associated Newton map $Nf$
 is indeterminate  at 0, and has a Cantor set of ``superstable" manifolds at the origin. 
 Notice that the treatment of the stable manifold theorem in \cite{ya} is covered 
 by remark \ref{rmk_Pesin}. 
  Let $\pi$ denote the blow up
at the origin, with exceptional divisor $E$. The lifted map $\widetilde{Nf}=
\pi^{-1}\rond Nf \rond \pi$
has two indeterminacy points  on $E$, whose 
image is $E$, and it is possible to reproduce the analysis 
of section \ref{sec_speed} --in our notations, $d_1=d_2=1$ and
 $\ell_1=\ell_2=2$-- hence the Cantor set of stable manifolds.

Using the formalism of horizontal-like maps and graph transform 
for currents allows for example to obtain parallel results in the case $P_2$ is 
arbitrary --this corresponds to arbitrary $d_i$'s. 
Of course the superstable manifolds have to be replaced by currents. \\

More generally, similar analysis should arise when  studying local dynamics
involving both some (super-)attracting behaviour and blow ups, e.g. near 
superattractive and (some) periodic indeterminacy points. 

\subsection{H{\'e}non(-like) maps.} 
Let $f$ be a regular polynomial automorphism
of $\cd$ (composition of complex H{\'e}non maps).
It is well known that there exists a 
bidisk $\bb=\dd_R^2$ such that points outside $\bb$ are wandering. The triple
$(f,\bb,\bb)$ is horizontal-like (H{\'e}non-like) in our sense. Assume 
$f(\bb)\cap\bb$ is a finite union of disjoint connected open sets
$U_1\cup\cdots\cup U_m$; if $L$ is a 
horizontal line, then $f(L)\cap\bb=L_1\cup\cdots\cup L_m$.
If all components $L_i$ have horizontal degree 1, $(f,\bb,\bb)$ is a complex 
horseshoe. In the general case each $U_i$ is equipped with a degree 
$d_i=\mathrm{deg}(L_i)$, independent of $L$. This discussion
is valid for H{\'e}non-like maps as well.\\

For instance this occurs when perturbing a polynomial in $\cc$ with disconnected 
Julia set: let for example
$p$ be a cubic polynomial in $\cc$, with escaping radius $R$, such that
one critical point $c_1$ satisfies $\abs{f(c_1)}>R$ and the other $c_2$ is attracted 
by a periodic sink.   Then $f^{-1}(\dd_R)$ has two components,
 one containing $c_2$ and the other containing no critical point, 
and for small $a$, 
the H{\'e}non map $f(z,w)=(aw+p(z), az)$ satisfies the 
previous   assumptions with $m=2$, $d_1=1$ and $d_2=2$. \\

As before it is possible to do symbolic dynamics by 
assigning to each point in $K=K^+\cap K^-$ (with the usual notation) its itinerary
in $\set{1,\ldots, m}^{\mathbb{Z}}$, and get a decomposition of $K$ into 
sets $K_\alpha$.

Now if $\alpha^+\in\set{1,\ldots, m}^{\nn^*}=\Sigma^+$ be a one sided 
symbol sequence, 
$$K_{\alpha^+}^+=\set{p\in\bb, ~\forall i\geq1,~ f^i(p)\in U_i}$$ is a vertical 
closed set carrying a current $T^+_{\alpha^+}$, and we get a decomposition 
$T^+=\int_{\Sigma^+} T^+_{\alpha^+} d\nu(\alpha^+)$, similarly to section 
\ref{sec_speed}. If some $d_i>1$, almost every $T^+_{\alpha^+}$ has bounded 
potential by proposition \ref{prop_apriori}.

If $\alpha\in\set{1,\ldots, m}^{\mathbb{Z}}$,  $K_\alpha=K^+_{\alpha^+}
\cap K^-_{\alpha^-}$ supports the probability measure $\mu_\alpha=
T^+_{\alpha^+}\wedge T^-_{\alpha^-}$, hence a decomposition of 
the maximal entropy measure $\mu=T^+\wedge T^-$.\\

It seems to be an interesting  question to describe in general the
extremal decomposition of the currents $T^{+/-}\rest{\bb}$. The latter examples 
provide a first insight into this question.

For another occurrence of sequences of H{\'e}non-like maps inside a given H{\'e}non map,
see \cite{bs10}.

\bigskip

\begin{small}
\noindent{\sc TC.D. and N.S. : Math{\'e}matique,
B{\^a}timent 425, 
Universit{\'e} de Paris Sud, 
91405 Orsay cedex, France.}\\
{\tt tiencuong.dinh@math.u-psud.fr, nessim.sibony@math.u-psud.fr}\\

\noindent{\sc R.D. : Institut de Math{\'e}matiques de Jussieu,
Universit{\'e} Denis Diderot,
Case 7012,
2 place Jussieu,
75251 Paris cedex 05, France.}\\
{\tt dujardin@math.jussieu.fr}\\
\end{small}

\end{document}